%% file: RealBinaryOctics-ArXiv.tex
\documentclass{amsart}

\usepackage{hyperref}
\usepackage{amssymb,amscd}
\usepackage{epic}
\usepackage[small,nohug,heads=littlevee]{diagrams}
\diagramstyle[labelstyle=\scriptstyle]

\title{On the Geometry of the Moduli Space of Real Binary Octics}
\author{Kenneth C. K. Chu}

\renewcommand*{\qed}{\mbox{} \hfill $\Box$}
\renewcommand*{\proof}{\noindent {\small P{\scriptsize ROOF}} \quad}
\renewcommand*{\i}{\mathbf{i}}
\usepackage{RealBinaryOctics}

\begin{document}

\maketitle
\begin{center} \vskip -0.8cm {\small \today} \end{center}

\begin{abstract}
\vskip 0.0cm
The moduli space of smooth real binary octics has five connected components.  They parametrize the real binary octics whose defining equations have $0, 1, \ldots, 4$ complex-conjugate pairs of roots respectively.  We show that the GIT-stable completion of each of these five components admits the structure of an arithmetic real hyperbolic orbifold.  The corresponding monodromy groups are, up to commensurability, discrete hyperbolic reflection groups, and their Vinberg diagrams are computed.  We conclude with a simple proof that the moduli space of GIT-stable real binary octics itself cannot be a real hyperbolic orbifold.
\end{abstract}



\input{RealBinaryOctics-ArXiv-body.tex}

\vskip 1.3cm
\noindent
{
\small
Department of Mathematics, University of Texas at Austin \\
1 University Station, C1200, Austin, Texas   78712, USA}  \\
{\tt chu@math.utexas.edu}


\end{document}

%% file: RealBinaryOctics-ArXiv-body.tex


\input{Introduction.tex}

\input{ComplexBallQuotient.tex}

\input{ACTconstruction.tex}

\input{CompareStabilizers.tex}

\input{InducedPermutation.tex}

\input{ComputationalResults.tex}

\input{NonOrbifoldPoints-nocurves.tex}

\input{Summary.tex}




\bibliographystyle{alpha}
\bibliography{PhDbibfile}


%% file: Introduction.tex
\section{Introduction}
\setcounter{theorem}{0}\setcounter{figure}{0}\setcounter{table}{0}

A \emph{(complex) binary octic} refers to a hypersurface of degree eight in the complex projective line $\CPone$.  One can think of a binary octic as an 8-point configuration in $\CPone$, counting multiplicity.  A binary octic is said to be \emph{smooth} if it is smooth as a hypersurface in $\CPone$; equivalently, it is smooth if the eight roots of any of its defining polynomials are pairwise distinct.  The GIT-stable (or more briefly, stable) binary octics are those with at worst triple-point singularities.  A \emph{real binary octic} is a binary octic that is preserved by complex conjugation on $\CPone$.

Using periods of certain branched covers of $\CPone$, Deligne-Mostow \cite{DM:hypergeom}, Terada \cite{Terada1985,Terada1983}, Matsumoto-Yoshida \cite{MY:complex8points} have described the arithmetic hyperbolic $5$-ball quotient structure of the moduli space $\Ms$ of stable complex binary octics.  The use of periods of curves is classical, for instance, in the construction of the moduli space of elliptic curves and Picard curves \cite{Picard1883}.  Kond$\bar{\textnormal{o}}$ \cite{Kondo:Ordered8Points} produced the same description of $\Ms$ using periods of $K$3 surfaces.

Following the approach of Allcock-Carlson-Toledo in \cite{ACT:realcubics} for real cubic surfaces and \cite{ACT:realsextics} for real binary sextics, this paper describes how the Deligne-Mostow construction of the moduli space of complex binary octics gives rise to an arithmetic real hyperbolic orbifold structure on the GIT-stable completion of each of the components of the moduli space of smooth real binary octics.  Unlike in \cite{ACT:realcubics} and \cite{ACT:realsextics}, the scalar ring involved here is the Gaussian integers and the lattice involved is no longer unimodular.  These lead to considerable added computational complexities, as well as the unforeseen semi-direct product structure of one of the monodromy groups.
Another key result in \cite{ACT:realcubics} (respectively \cite{ACT:realsextics}) is that the completions of the individual components of the moduli space of smooth real cubic surfaces (respectively real binary sextics) glue together nicely to give the moduli space of GIT-stable real cubic surfaces (respectively real binary sextics) the structure of a \emph{non-arithmetic} real hyperbolic orbifold.  This nice property is however not shared by 
the moduli space of stable real binary octics.  In fact, this paper concludes with a simple proof that the latter space cannot be a real hyperbolic orbifold at all.




%% file: ComplexBallQuotient.tex
\section{The Moduli Space of Complex Binary Octics as an Arithmetic Quotient of $\CHfive$ }
\label{section:SpaceOfFramedForms}
\setcounter{theorem}{0}\setcounter{figure}{0}\setcounter{table}{0}


Let $\Po$ be the space of smooth binary octic forms (homogeneous binary polynomials of degree eight) and $\Ps$ be the space of stable binary octic forms.  We take the covering $\Fo\rightarrow\Po$ corresponding to the kernel $\P\Gamma$ of a certain representation $\pi_{1}(\Po)\overset{\rho}{\longrightarrow}\Isom(\CHfive)$.  Then, $\P\Gamma$ acts on $\Fo$ as deck transformations, and on $\CHfive$ via $\rho$.  On the other hand, a certain group $G$, derived from $\PGLtwoC$, acts on $\Ps$, preserving $\Po$, such that $\Ms \cong \Ps/G$ and $\Mo \cong \Po/G$.   Let $\Fs\rightarrow\Ps$ be the Fox completion \cite{Fox} of $\Fo\rightarrow\Po$.  Then, the actions of $G$ and $\P\Gamma$ on $\Fo$ naturally extend to $\Fs$ such that $\Fs/G \cong \CHfive$ and $\P\Gamma\backslash\Fs\cong\Ps$.  Roughly, the complex $5$-ball quotient structure on $\Ms$ arises as follows:  $\Ms$ $\cong$ $\Ps/G$ $\cong$ $\left(\P\Gamma\backslash\Fs\right)/G$ $\cong$ $\P\Gamma\backslash\left(\Fs/G\right)$ $\cong$ $\P\Gamma\backslash\CHfive$.  We remark that, in this way, $\Ms$ and $\P\Gamma\backslash\CHfive$ are isomorphic as complex analytic (quasi-projective) varieties, but not as complex analytic orbifolds. 

In this section, we give some details of the above constructions and state the known properties of the ball quotient $\P\Gamma\backslash\CHfive$ which will be useful in the sequel.  We refer to the literature for proofs whenever possible.


\subsection{The fibration of cyclic covers branched over octics and the Hermitian structure of the cohomology of its fiber}
\setcounter{theorem}{0}\setcounter{figure}{0}\setcounter{table}{0}

Define
\begin{equation*}
   \scriptX \; := \;
   \left\{\;
      (\; p \,,\,[x_{0}:x_{1}:y] \;) \in \mathcal{P}\times\P(1,1,2)
      \;\;\left|\;\;
      y^{4} - p(x_{0},x_{1}) = 0
      \;\right.
   \right\},
\end{equation*}
where $\mathcal{P}$ is the space of all binary octic forms and $\P(1,1,2)$ is the weighted projective space of weights $(1,1,2)$.  Let $\i$ denote $\sqrt{-1}$.  Define the maps
\begin{equation*}
\begin{array}{ccccccl}
   \sigma : & \scriptX & \longrightarrow & \scriptX : &
              (\,p\,,\,[x_{0}:x_{1}:y]\,) & \longmapsto &
                \left(\,p\,,\,[x_{0}:x_{1}:\i\,y]\,\right), \\ \\
   \Pi :    & \scriptX & \longrightarrow  & \mathcal{P} : &
              (\,p\,,\,[x_{0}:x_{1}:y]\,) & \longmapsto  & p, \\ \\
   \pi :    & \scriptX & \longrightarrow  & \CPone : &
              (\,p\,,\,[x_{0}:x_{1}:y]\,) & \longmapsto  & [x_{0}:x_{1}].
\end{array}
\end{equation*}
Let $\scriptXo := \Pi^{-1}(\Po)$.  Then, for each $p \in \Po$, the fiber
\begin{equation*}
     X_{p} := \Pi^{-1}(p)
   = \left\{\; [x_{0}:x_{1}:y] \in \P(1,1,2) \;|\; y^{4} - p(x_{0},x_{1}) = 0 \;\right\}
\end{equation*}
is a (smooth) compact Riemann surface.  The map $\sigma : \scriptX \longrightarrow \scriptX$ is a cyclic
action on $\scriptX$ of order 4.  $\sigma$ preserves every fiber of $\Pi$, hence restricting to a cyclic action of order 4 on each fiber $X_{p} := \Pi^{-1}(p)$, $p \in \Po$.  The map $\pi : \scriptX \longrightarrow \CPone$ is well-defined since $[0:0:1] \in \P(1,1,2) - \scriptX$.  Observe that for each $p \in \Po$, the restricted map $\pi|_{X_{p}} : X_{p} \longrightarrow \CPone$ is a cyclic cover of $\CPone$ of degree 4 branched over the eight distinct roots of $p(x_{0},x_{1})$ in $\CPone$, and it has exactly eight ramification points, each with ramification index 4.  By the Riemann-Hurwitz theorem, $g(X_{p}) = h^{1,0}(X_{p})  = 9$, for each $p \in \Po$.  Thus, $\scriptXo \overset{\Pi}{\longrightarrow} \Po$ is a fibration whose fiber over each $p \in \Po$ is the compact Riemann surface $X_{p} := \Pi^{-1}(p)$, which has genus 9 and is a cyclic covering of $\CPone$ branched over the roots in $\CPone$ of the polynomial $p(x_{0},x_{1})$.

Next, for each $p \in \Po$, define
\begin{equation*}
   \LXp \; := \; H^{1}_{\sigma^{2}=-1}(X_{p},\Z) \; := \;
   \left\{\;
      \phi \in H^{1}(X_{p},\Z) \;\;|\;\; \sigma^{2}(\phi) = - \phi
   \;\right\}.
\end{equation*}
Then $\sigma|_{\LXp}$ satisfies $\sigma^{2} + 1 = 0$.  Consequently, if we define multiplication by $-\i$ in
$\LXp$ by
\begin{equation*}
   -\i\cdot\phi := \sigma(\phi),
\end{equation*}
then $\LXp$ becomes a $\Zi$-module.

\begin{proposition}\label{GaussianModule}\quad
   With the above $\Zi$-module structure, $\LXp$  
   becomes a free $\Zi$-module of rank 6.
\end{proposition}
\proofoutline  Torsionfree-ness
of $\LXp$ over $\Z$ readily implies its torsionfree-ness
over $\Zi$.  Since $\Zi$ is a PID, $\LXp$ is a free $\Zi$-module.
By examining the complex-valued de Rham cohomology of
$X_{p}$, we find that $\rank_{\Zi}\left(\LXp\right) = 6$. \qed

\vskip 0.5cm
\noindent
Consider the embedding $\LXp \hookrightarrow H^{1}_{\sigma=-\i}(X_{p},\C)$
induced by
\begin{diagram}
H^{1}_{\sigma^{2}=-1}(X_{p},\Z)&\rInto^{\mbox{}\;\quad\quad\mbox{}}&H^{1}_{\sigma^{2}=-1}(X_{p},\Z)\otimes_{\Z}\C&\rTo^{\mbox{}\quad\sim\!\!\!\!\mbox{}}&H^{1}_{\sigma^{2}=-1}(X_{p},\C) \\
\parallel & & & & \parallel \\ 
\LXp & & & & \textnormal{\small$H^{1}_{\sigma=-\i}(X_{p},\C) \oplus H^{1}_{\sigma=\i}(X_{p},\C)$} \\
        & \rdInto(4,2)& & & \dOnto \\
        & & & & H^{1}_{\sigma=-\i}(X_{p},\C)  
\end{diagram}
Let
$h' : H^{1}_{\sigma=-\i}(X_{p},\C) \times H^{1}_{\sigma=-\i}(X_{p},\C) \longrightarrow \C$
be the Hermitian form given by
\begin{equation*}
   (\alpha,\beta) \overset{h'}{\longmapsto} 2\,\i\,\int_{X_p}\, \alpha\wedge\overline{\beta}.
\end{equation*}
The above Hermitian form induces a Gaussian lattice structure on $\LXp$, as the following Proposition shows:

\begin{proposition}\quad
\label{Proposition:Lorentzian}
\begin{enumerate}

\item\label{Lorentzian}
$\dim_{\C}H^{1,0}_{\sigma=-\i}(X_{p},\C) = 1$, and
$\dim_{\C}H^{0,1}_{\sigma=-\i}(X_{p},\C) = 5$.
$h'$ is positive-definite on $H^{1,0}_{\sigma=-\i}(X_{p},\C)$
and negative-definite $H^{0,1}_{\sigma=-\i}(X_{p},\C)$. 
Consequently,
$\left(\;
      H^{1}_{\sigma=-\i}(X_{p},\C) \,,\, h'
   \;\right)$
is isometric to the standard Lorentzian-Hermitian space
$\C^{1,5}$ $=$ $\C^{1+,5-}$. 

\item\label{LatticeStructureOfLXp}
  Let $h$ be the pull-back to $\LXp$ of the Lorentzian-Hermitian form
  $h' : H^{1}_{\sigma=-\i}(X_{p},\C) \times H^{1}_{\sigma=-\i}(X_{p},\C) \longrightarrow \C$
  by the embedding $\LXp \hookrightarrow H^{1}_{\sigma=-\i}(X_{p},\C)$.  Then,
  $h$ is in fact $\Zi$-valued on $\LXp \times \LXp$, and it is a
  $\Zi$-Hermitian form on $\LXp$ given by the following formula:
  \begin{equation*}
     h(\,\xi\,,\,\eta\,) \; = \; -\,\Omega(\,\xi \,,\,\sigma(\eta)\,) - \i\,\Omega(\,\xi\,,\,\eta\,),
     \quad\textnormal{for any} \;\, \xi, \; \eta \in \LXp,
  \end{equation*}
  where $\Omega : H^{1}(X_{p},\Z) \times H^{1}(X_{p},\Z) \longrightarrow \Z$
  is given by
  \begin{equation*}
     \Omega(\,\alpha\,,\,\beta\,) := \langle\; \alpha \cup \beta \;,\; [X_{p}]\;\rangle.
  \end{equation*}

\item\label{explicitLXp}
  The Lorentzian $\Zi$-Hermitian quadratic form on $\LXp$
  constructed as in (\ref{LatticeStructureOfLXp})
  is abstractly isometric to the following $\Zi$-lattice:
  \begin{equation*}
     \Lambda \; := \;
     \left(\;\Zi^{6}\;,\;
     \textnormal{\small$\cuspsummand\oplus\cuspsummand\oplus\Htwosummand$}
     \;\right).
  \end{equation*}
\end{enumerate}  
\end{proposition}
\proofremark  The proofs of all three statements follow by direct computations.  The proof of \eqref{LatticeStructureOfLXp} is similar to that in paragraph (4.3) in \cite{ACT:complexcubics}, while that of \eqref{explicitLXp} can be easily inferred from the results in \cite{MY:complex8points}.\qed


\subsection{The space of framed octic forms}
\setcounter{theorem}{0}\setcounter{figure}{0}\setcounter{table}{0}

In this section, we describe the space of framed smooth octic forms and its Fox completion \cite{Fox}, the space of framed stable octic forms.  They are the domains of the period maps described in the subsequent sections.  The complex ball quotient structure of $\Ms$ arises through these period maps.  We omit all proofs, but refer to \cite{ACT:complexcubics}, which treats the analogous case of the complex cubic surfaces.

\begin{definition}\label{defn:framedsmoothform}
  A \emph{framed smooth octic form} over $p \in \Po$ is a
  ``projective equivalence class'' of an (abstract)
  isometry of $\LXp \overset{\sim}{\longrightarrow} \Lambda$,
  where two such isometries are said to be
  ``projectively equivalent'' if one is a $\Zi$-unit scalar
  multiple of the other.
\end{definition}

Let $\Lambda(\scriptXo)$ be the sheaf over $\Po$ associated to the presheaf $U \mapsto H^{1}_{\sigma^{2}=-1}(\Pi^{-1}(U),\Z)$.  Proposition \ref{Proposition:Lorentzian}\eqref{explicitLXp} implies that $\Lambda(\scriptXo)$ is a sheaf over $\Po$ of $\Zi$-Hermitian modules, with stalks isomorphic to the rank-six $\Zi$-lattice $\Lambda$.  Let $\P\textnormal{Hom}(\Lambda(\scriptXo),\Po\times\Lambda)$ be the sheaf of projective equivalence classes of sheaf homomorphisms from $\Lambda(\scriptXo)$ to $\Po\times\Lambda$.

\begin{definition}\quad  The \emph{space $\Fo$ of framed smooth octic forms} over $\Po$ is the subsheaf of \\
$\nobreak{\P\textnormal{Hom}(\Lambda(\scriptXo),\Po\times\Lambda)}$ consisting of projective equivalence classes of sheaf homomorphisms $\Lambda(\scriptXo)\rightarrow\Po\times\Lambda$ which restrict to an isometry on each stalk.
\end{definition}

$\Fo$ is a complex manifold and its stalks are the framed smooth octic forms, as defined in Definition
\ref{defn:framedsmoothform}.   $\Fo$ can be alternatively described as the Galois covering of $\Po$ associated to the kernel of the ``projectivized monodromy representation''
\begin{equation*}
   \P\rho : \pi_{1}(\Po,p_{0}) \longrightarrow \P\Isom(\Lambda(X_{p_{0}})) \cong \P\Isom(\Lambda),
\end{equation*}
which of course derives from the standard monodromy representation
\begin{equation*}
    \rho : \pi_{1}(\Po,p_{0}) \longrightarrow \Isom(\Lambda(X_{p_{0}})),
\end{equation*}    
where $p_{0}\in\Po$ is an arbitrary but fixed smooth octic.  It is clear from this description of $\Fo$ as a Galois covering over a path-connected base space that it is connected.  The monodromy group --- and the deck transformation group --- $\rho(\pi_{1}(\Po,p_{0}))$ $\subset$ $\P\Isom(\Lambda)$ turns out to be all of $\P\Isom(\Lambda)$.  So, $\P\Gamma := \P\Isom(\Lambda)$ acts on $\Fo$ as deck transformations, and $\P\Gamma\backslash\Fo \cong \Po$.

Let $G := \GL(2,\C) / \langle \textnormal{all eighth roots of unity} \rangle$.  $G$ acts naturally on $\Po$ (by ``linear change of variables") and this action extends to a free action on $\Fo$ via induced action on cohomology.

Next, let $\Ps$ be the space of all stable binary octic forms and $\Fs$ be the Fox completion (see \cite{Fox}) of the covering $\Fo\rightarrow\Po$.  $\Fs$ is a branched covering of $\Ps$ with four-fold branching over $\Delta^{1}_{s} \subset \Ps$, the locus in $\Ps$ corresponding to octics with one double
point and no other singularities.  Intuitively, $\Fs$ coincides with $\Fo$ over $\Po$, and, for a singular octic $p\in\Delta^{1}_{s}$, $\Fs$ retains information about the vanishing cohomology corresponding to the singularities of $p$.  We call $\Fs$ the \emph{space of framed stable octic forms}.

The actions of $G$ and $\P\Gamma$ on $\Fo$ extend naturally to $\Fs$, and it can be shown that $\P\Gamma\backslash\Fs \cong \Ps$.


\subsection{The complex period map and the $\CHfive$ quotient structure of $\Ms$}
\label{section:complexperiodmap}
\setcounter{theorem}{0}\setcounter{figure}{0}\setcounter{table}{0}

The period map of interest to us is defined as follows:
\begin{equation*}
\begin{array}{rcl}
  \Fo & \overset{\p}{\longrightarrow} & \CHfive = \CH(\Lambda\otimes_{\Zi}\C) \\
  \left[\LXp\overset{i}{\rightarrow}\Lambda\right] & \longmapsto     & i(H^{1,0}_{\sigma=-\i}(X_{p}))
\end{array}.
\end{equation*}

Note that $\P\Gamma = \P\Isom(\Lambda)$ naturally acts on $\CHfive = \CH(\Lambda\otimes_{\Zi}\C)$.  The period map $\p$ turns out to be holomorphic, invariant under the action of $G$ on $\Fo$, and it is equivariant with respect to the actions of $\P\Gamma = \P\Isom(\Lambda)$ on $\Fo$ and $\CHfive$.

The period map $\p$ extends holomorphically to $\Fs$ to a $(G\acts\Fs)$-invariant and $\P\Gamma$-equivariant map, also denoted by $\p$.  The map $\p$ therefore descends to a map $\p : \Fs/G \longrightarrow \CHfive$, which turns out to be an isomorphism of complex manifolds.  Furthermore, $\p$ maps $\Fo$ bijectively to $(\CHfive - \mathcal{H})$, where
\begin{equation*}
   \mathcal{H} := \bigcup\;\left\{\; \CH(r^{\perp})\subset\CHfive
                  \;\left\vert\;
                  \begin{array}{c}
                     \textnormal{\small $r$ is a vector in $\Lambda$} \\
                     \textnormal{\small of squared norm $-2$}
                  \end{array}
                  \right.\;\right\},
\end{equation*}
restricting also to an isomorphism of complex manifolds $\Fo/G \overset{\p}{\longrightarrow} (\CHfive-\mathcal{H})$.

The results of Deligne-Mostow \cite{DM:hypergeom} and Matsumoto-Yoshida \cite{MY:complex8points} show that $\Ms$ and $\P\Gamma\backslash\CHfive$ are isomorphic as complex analytic (quasi-projective) varieties via the following series of isomorphisms:
\begin{equation*}
  \Ms :=    \P(\Ps)/\P\GL(2,\C) \cong \Ps/G
      \cong (\P\Gamma\backslash\Fs)/G
      \cong \P\Gamma\backslash(\Fs/G) 
      \cong \P\Gamma\backslash\CHfive.
\end{equation*}
We stress that $\Ms$ and $\P\Gamma\backslash\CHfive$
are isomorphic only as complex analytic varieties, but
not as complex analytic orbifolds.

%% file: ACTconstruction.tex
\section{The Allcock-Carlson-Toledo Construction of $\MsR$}
\label{section:BasicFactsOfRealBinaryOctics}
\setcounter{theorem}{0}\setcounter{figure}{0}\setcounter{table}{0}

As shown in the last section, the moduli space $\Ms$ of stable binary octics is isomorphic as complex analytic varieties to the ball quotient $\P\Gamma\backslash\CHfive$.  It turns out that periods in $\CHfive$ corresponding to real octics lie on a certain collection of copies of real hyperbolic $5$-space $\RHfive$ inside $\CHfive$.  Roughly speaking, the Allcock-Carlson-Toledo construction of $\MsR$ is simply to extract this collection of copies of $\RHfive$ and re-assemble them according to the ``expected'' quotient structure of $\MsR$.


\subsection{Complex conjugation and the antipodal map on $\CPone$ and their related maps}
\setcounter{theorem}{0}\setcounter{figure}{0}\setcounter{table}{0}

We are interested in the moduli space of stable real binary octics, namely, binary octics whose coefficients are real.  We can also view these as those (a priori complex) octics that are invariant under the action induced on the space of binary octic forms by the ``usual'' antiholomorphic involution on $\CPone$, i.e. complex conjugation $\kappa : \CPone \longrightarrow \CPone$.

Up to $\P\GLtwoC$-conjugacy, there is exactly one more antiholomorphic involution on $\CPone$, namely, the antipodal map $\alpha : \CPone \longrightarrow \CPone$.  (See \cite{Kollar:realforms}.)  For reasons that will become apparent shortly, we need to deal with the octics which are preserved by $\alpha$ as well.

\begin{definition}\label{kappaalpha}\quad
Define the maps $\kappa:\C^{2}\longrightarrow\C^{2}$, and $\alpha:\C^{2}\longrightarrow\C^{2}$ respectively by\;\; $\kappa(x_{0},x_{1}) \; := \; (\; \overline{x_{0}} \;, \overline{x_{1}}\;)$, \;\;and\;\; $\alpha(x_{0},x_{1}) \; := \; (\; \overline{x_{1}} \;,-\overline{x_{0}}\;)$.
\end{definition}

\begin{definition}\label{ActionOnOcticForms}\quad
Let $\nu:\C^{2}\longrightarrow\C^{2}$ be either $\kappa$ or $\alpha$ as in Definition \ref{kappaalpha}.
We define the action of $\nu$ on the space of complex binary octic forms $\mathcal{P}$ as follows:
\begin{equation*}
   (\nu\cdot p)(x_{0},x_{1}) \;\; := \;\; \overline{p(\nu(x_{0},x_{1}))},
   \quad \textnormal{for}\;\; p \in \mathcal{P}.
\end{equation*}
\end{definition}

\begin{remark}\quad
$\kappa$ descends to complex conjugation on $\CPone$, whereas $\alpha$ descends to the antipodal map on $\CPone$.
We will also use $\kappa$ to denote complex conjugation on $\CPone$ and $\alpha$ the antipodal map on $\CPone$.  Which map is intended should be clear from the context.
\end{remark}

\begin{definition}\quad  A binary octic form is said to be \emph{real} (respectively \emph{antipodal}) if it is preserved by complex conjugation $\C^{2}\overset{\kappa}{\longrightarrow}\C^{2}$ (respectively the antipodal map $\C^{2}\overset{\alpha}{\longrightarrow}\C^{2}$) via the action as in Definition \ref{ActionOnOcticForms}. We denote by $\PoR$ the set of smooth real binary octic forms, and by $\Poa$ the set of smooth antipodal binary octic forms.   We denote by $\FoR$ and $\Foa$ the preimages of $\PoR$ and $\Poa$, respectively, under the covering map $\Fo \longrightarrow \Po$.
\end{definition}

\begin{remark}\quad
There are smooth octics that are preserved by both complex conjugation and the antipodal map.  In other words, $\PoR \cap \Poa \neq \varnothing$.
\end{remark}

\begin{definition}\quad
Let $\GL(2,\C)'$ be the group of all linear and antilinear automorphisms of $\C^{2}$; note that $\GLtwoC' = \GLtwoC\rtimes\langle\kappa\rangle$.  Let every linear element $g \in \GL(2,\C)'$ and every antilinear element $\nu \in \GL(2,\C)'$ act on $\C^{3}$ respectively by:
\begin{equation*}
  g(x_{0},x_{1},y) \; := \; \left(\, g(x_{0},x_{1}) \,,\, y \,\right),
  \quad\textnormal{and}\quad
  \nu(x_{0},x_{1},y) \;:= \; \left(\, \nu(x_{0},x_{1}) \,,\, \overline{y} \,\right),
\end{equation*}
We will also consider elements of $\GLtwoC'$ as automorphisms of $\P(1,1,2)$ via the representation $\GLtwoC'\longrightarrow\Aut'\,\P(1,1,2)$ corresponding to the action $\GLtwoC'\act\C^{3}$ above, where $\Aut'\,\P(1,1,2)$ is the automorphism group of $\P(1,1,2)$ induced by linear and antilinear automorphisms of $\C^{3}$.
\end{definition}

\begin{definition}\quad
Let $\GR$ be the centralizer $\mathcal{C}_{\Aut\,\P(1,1,2)}(\kappa)$ of $\kappa \in \Aut'\,\P(1,1,2)$ in $\Aut\,\P(1,1,2) \subset \Aut'\,\P(1,1,2)$.  Let $\Ga$ be the centralizer $\mathcal{C}_{\Aut\,\P(1,1,2)}(\alpha)$ of $\alpha\in\Aut'\,\P(1,1,2)$ in $\Aut\,\P(1,1,2) \subset \Aut'\,\P(1,1,2)$.
\end{definition}

Straightforward calculations show that $\GR = \GLtwoR/\langle\pm 1\rangle$ and
\begin{eqnarray*}
     \Ga &=& \left\{\;g\in\GLtwoC
                   \;\,\left\vert\;\;
                   {\textnormal{\scriptsize$\left[\begin{array}{rr}0&1\\ \!\!\!-1&0\end{array}\right]$}}
                   \cdot\overline{g}\,=\,\pm\,g\cdot
                   {\textnormal{\scriptsize$\left[\begin{array}{rr}0&1\\ \!\!\!-1&0\end{array}\right]$}}
                   \,\right.\right\} \\
            &=& \left\{\;\left.
                   {\textnormal{\scriptsize$\left[\begin{array}{rr}z_{1}&z_{2} \\ \!\!\!\pm\,\overline{z_{2}}&\mp\,\overline{z_{1}}\end{array}\right]$}}
                   \in\C^{2 \times 2}\;\,\right\vert\;\,
                   |z_{1}|^{2} + |z_{2}|^{2} \neq 0
                   \,\right\}.
\end{eqnarray*}

By an \emph{anti-isometry} on a $\Zi$-Hermitian lattice $(\,V\,,\,\langle\cdot,\cdot\rangle\,)$ (or a complex vector space equipped with a Hermitian inner product), we mean a bijective antilinear map $\nu : V \longrightarrow V$ such that $\langle \nu(x), \nu(y) \rangle = \overline{\langle x, y \rangle}$, for all $x, y \in V$.

\begin{definition}\quad
Let $\Fo'$ be the space of all pairs $(p,[i])$, where $p \in \Po$, $\LXp\overset{i}{\longrightarrow}\Lambda$ is either an isometry or an anti-isometry, and $[i]$ is the projective equivalence class of $i$.  Let every linear element $g \in \GLtwoC'$ and every antilinear element $\nu\in\GLtwoC'$ act on $\Fo'$ respectively by
\begin{equation*}
  (\,p\,,\,[i]\,)\cdot g \; := \; \left(\, p\circ g\,,\,[i\circ(g^{*})^{-1}]\,\right),
  \quad\textnormal{and}\quad
  (\,p\,,\,[i]\,)\cdot \nu \; := \; \left(\, \overline{p\circ \nu}\,,\,[i\circ(h^{*})^{-1}]\,\right).
\end{equation*}
\end{definition}

Note that, for $p\in\PoR$ and $g\in\GR$ (respectively $p\in\Poa$ and $g\in\Ga$), we have the the following commutative diagrams:
\begin{equation*}
\begin{CD}
     X_{p \circ g} @> \kappa_{p \circ g} >> X_{p \circ g} \\
     @V g VV                                              @VV g V       \\
     X_{p}           @>> \kappa_{p}           >  X_{p}           \\
\end{CD}
\quad\quad\quad\quad\quad\quad\quad\quad
\begin{CD}
     X_{p \circ g} @> \alpha_{p \circ g} >> X_{p \circ g} \\
     @V g VV                                              @VV g V       \\
     X_{p}           @>> \alpha_{p}           >  X_{p}           \\
\end{CD}
\end{equation*}


\subsection{The deformation types of real and antipodal smooth octics and forms}
\label{section:TopologicalTypes}
\setcounter{theorem}{0}\setcounter{figure}{0}\setcounter{table}{0}

There are five distinct deformation types of smooth real binary octics, in the sense that a real octic, of any fixed deformation type, cannot be deformed to a real octic of a different type through the space $\OoR = \left.\PoR\right/\Re^{*}$ of smooth real octics (where $\Re^{*} := \Re\,\backslash\{0\}$ acts by scalar multiplication on the real octic forms, as usual).  In other words, $\OoR$ has five connected components, i.e.
\begin{equation*}
  \OoR \;\; = \;\; 
  \OoRzero \,\bigsqcup\, \OoRone \,\bigsqcup\, \OoRtwo \,\bigsqcup\, \OoRthree \,\bigsqcup\, \OoRfour,
\end{equation*}
where $\OoRzero, \ldots, \OoRfour$ parametrize the five types of real binary octics according Table \ref{table:DeformationTypesOfRealBinaryOctics}.
\begin{table}
\begin{center}
\begin{tabular}{|c||c|c|c|c|c|}
  \hline
  components of $\OoR$         & $\OoRzero$ & $\OoRone$ & $\OoRtwo$ & $\OoRthree$ & $\OoRfour$ \\ \hline
  {\# complex conjugate pairs} & 0 & 1 & 2 & 3 & 4 \\ \hline
  {\# real points}             & 8 & 6 & 4 & 2 & 0 \\ \hline
\end{tabular}
\vskip 0.1cm
\caption{Deformation types smooth real binary octics}
\label{table:DeformationTypesOfRealBinaryOctics}
\end{center}
\end{table}

On the other hand, every smooth antipodal octic can be deformed to every other smooth antipodal octic through smooth antipodal octics.  In other words, $\Ooa$ is connected and there is only one deformation type of smooth antipodal octics.

\begin{definition}\quad
Let $\MoR$ be the moduli space of smooth real binary octics and $\MoRzero, \MoRone, \ldots \MoRfour$ its five connected components of $\MoR$, parametrizing octics in $\OoRzero, \OoRone, \ldots, \OoRfour$, respectively.  (\,Therefore, $\MoR \;\; = \;\; \bigsqcup_{i=0}^{4}\MoRi$.\,)  Let $\Moa$ be the moduli space of smooth antipodal octics.
\end{definition}

By contrast, in order to count the number of connected components of $\PoR$, we need to take into account the fact that $\Re^{*}$ has two connected components.  Write $\PoRi$ for the preimage of $\OoRi$ under the projection $\PoR\longrightarrow\OoR=\left.\PoR\right/\Re^{*}$, $i = 0, \ldots, 4$.  Consider a smooth real binary octic in $\OoRi$, determined by say the roots of an octic form $p(x_{0},x_{1})\in\PoRi$.  Then, both $p(x_{0},x_{1})$ and $-p(x_{0},x_{1})$ descend to the same given octic ($8$-point configuration), but they may or may not belong to the same connected component of $\PoRi$.  It is now clear that each $\PoRi$, $i=0,\ldots,4$, has either one or two connected components, depending on whether or not any (hence every) element $p(x_{0},x_{1}) \in \PoRi$  can be deformed to its negative $-p(x_{0},x_{1})$ within $\PoRi$.  We now prove:

\begin{lemma}\label{DeformationTypesOfForms}\quad
$\PoRfour$ has two connected components\footnote{The author wishes to express his gratitude to Dr. J\'{a}nos Koll\'{a}r for pointing out the author's earlier overlooking of this fact in a private communication.}, whereas each of $\PoRzero$, $\PoRone$, $\PoRtwo$, $\PoRthree$, and $\Poa$ is connected.
\end{lemma}
\proof If we regard $x_{0}$ and $x_{1}$ as real variables, then each pair $p(x_{0},x_{1}), -p(x_{0},x_{1}) \in \PoRfour$ can be regarded as continuous $\Re$-valued nowhere vanishing functions of the real variables $x_{0}$, $x_{1}$ of opposite signs.  Consequently, any continuous deformation from $p(x_{0},x_{1})$ to $-p(x_{0},x_{1})$ through the space of continuous $\Re$-valued functions must pass through one that admits zeroes, thereby passing outside $\PoRfour$, since every smooth real binary octic form in $\PoRfour$ has no real roots.  This proves that $\PoRfour$ has two connected components.

Next, consider the following $1$-parameter family of binary polynomials:
\begin{equation*}
q_{3}(x_{0},x_{1};\theta) := (x_{0}\cos\theta-x_{1}\sin\theta) (x_{0}\sin\theta+x_{1}\cos\theta),
\quad \theta \in [0,\pi/2].
\end{equation*}
Then, $q_{3}(x_{0},x_{1};0)$ $=$ $x_{0}x_{1}$, whereas $q_{3}(x_{0},x_{1};\pi/2) = - x_{0}x_{1}$.  Let $r(x_{0},x_{1})$ be any smooth real binary sextic form with no real roots.  Then, $p(x_{0},x_{1};\theta) := q_{3}(x_{0},x_{1};\theta) r(x_{0},x_{1})$, $\theta \in [0,\pi/2]$, is a continuous path in $\PoRthree$ such that $p(x_{0},x_{1};0) = x_{0}x_{1}\!\cdot\!r(x_{0},x_{1})$, while $p(x_{0},x_{1};\pi/2)=-x_{0}x_{1}\!\cdot\!r(x_{0},x_{1})$.  This proves that $\PoRthree$ is connected.

Similarly, we may define continuous paths in $\PoRi$, $i=0,1,2$, whose endpoints are negatives of each other by using the following three families in place of $q_{3}$:
\begin{equation*}
\begin{array}{l}
q_{2}(x_{0},x_{1};\theta_{2}) \; := \\
\left(x_{0}\cos\theta_{2} - x_{1}\sin\theta_{2}\right)\left(x_{0}\sin\theta_{2} + x_{1}\cos\theta_{2}\right) \\
\times\left(x_{0}\cos(\theta_{2}+\pi/4) - x_{1}\sin(\theta_{2}+\pi/4)\right)\left(x_{0}\sin(\theta_{2}+\pi/4) + x_{1}\cos(\theta_{2}+\pi/4)\right),
\end{array}
\end{equation*}
\begin{equation*}
\begin{array}{l}
q_{1}(x_{0},x_{1};\theta_{1}) \; := \\
\overset{2}{\underset{n=0}{\textnormal{\Large$\prod$}}}
\left(x_{0}\cos(\theta_{1}+n\pi/6) - x_{1}\sin(\theta_{1}+n\pi/6)\right)
\left(x_{0}\sin(\theta_{1}+n\pi/6) + x_{1}\cos(\theta_{1}+n\pi/6)\right),
\end{array}
\end{equation*}
\vskip 0cm
\begin{equation*}
\begin{array}{l}
q_{0}(x_{0},x_{1};\theta_{0}) \; := \\
\overset{3}{\underset{n=0}{\textnormal{\Large$\prod$}}}
\left(x_{0}\cos(\theta_{0}+n\pi/8) - x_{1}\sin(\theta_{0}+n\pi/8)\right)
\left(x_{0}\sin(\theta_{0}+n\pi/8) + x_{1}\cos(\theta_{0}+n\pi/8)\right),
\end{array}
\end{equation*}
where $\theta_{2} \in [0,\pi/4]$, $\theta_{1}\in[0,\pi/6]$, $\theta_{0}\in[0,\pi/8]$.  Thus, $\PoRzero$, $\PoRone$, and $\PoRtwo$ are connected.  Lastly, we conclude that $\Poa$ is also connected by noting that $q_{0}(x_{0},x_{1};\theta_{0})$ is a family of antipodal octic forms (in addition to being real).
\qed

\vskip 0.5cm
In summary, $\PoR$ has six connected components, i.e.,
\begin{equation*}
  \PoR \;\; = \;\; 
  \PoRzero \,\bigsqcup\, \PoRone \,\bigsqcup\, \PoRtwo \,\bigsqcup\, \PoRthree \,\bigsqcup\,
  \PoRfourplus \,\bigsqcup\, \PoRfourminus,
\end{equation*}
where $\PoRfourplus$ and $\PoRfourminus$ are the two connected components of $\PoRfour$.



\subsection{Each $p\in\PoR\bigsqcup\Poa$ gives rise to an involutive anti-isometry of $\LXp$}
\setcounter{theorem}{0}\setcounter{figure}{0}\setcounter{table}{0}

If $p \in \PoR$, then complex conjugation $\CPone\overset{\kappa}{\longrightarrow}\CPone$ induces an antiholomorphic involution $\kappa_{p}$ on $X_{p}$.  Similarly, if $p \in \Poa$, then the antipodal map $\CPone\overset{\alpha}{\longrightarrow}\CPone$ likewise induces an antiholomorphic involution $\alpha_{p}$ on $X_{p}$.
For each octic $p \in \PoR\cap\Poa$, both $\kappa_{p}$ and $\alpha_{p}$ on $X_{p}$ are defined.

Let $p\in\PoR\cup\Poa$, and let $\nu_{p}$ be $\kappa_{p}$ or $\alpha_{p}$, whichever is defined on $X_{p}$.  Then the antiholomorphic involution $X_{p}\overset{\nu_{p}}{\longrightarrow}X_{p}$ induces an antilinear involution on $H^{1}(X_{p},\C)$ via 
\begin{equation*}
\begin{array}{ccl}
  H^{1}(X_{p},\C) & \overset{\nu'_{p}}{\longrightarrow} & H^{1}(X_{p},\C)               \\
  \phi            & \longmapsto                         & \overline{(\nu_{p})^{*}(\phi)}
\end{array}
\end{equation*}

\begin{lemma} \quad
\begin{enumerate}
    \item   The map $\nu'_{p}$ preserves both the Hodge decomposition and the
               $\sigma$-eigenspace decomposition of $H^{1}(X_{p},\C)$.
    \item   The antiholomorphic map $\nu'_{p}$ restricts to an involutive
               anti-isometry on $H^{1}_{\sigma=-\i}(X_{p},\C)$, which in turn restricts
               to an involutive anti-isometry on the $\Zi$-lattice on $\LXp$.
\end{enumerate}
\end{lemma}
\proofoutline Since $\nu_{p}$ is antiholomorphic, the pullback $\nu_{p}^{*}$ switches Hodge types of $\C$-valued differential forms; similarly, complex conjugation on $\C$-valued differential forms switches Hodge types.  Hence, $\nu'_{p}$ preserves Hodge types.  To prove that $\nu'_{p}$ preserves $\sigma$-eigenspaces, we first state two facts: $\sigma \circ \nu_{p} = \nu_{p} \circ \sigma^{3}$, and that the action of
$\sigma^{*}$ on $\C$-valued differential forms commutes with complex conjugation of differential forms.  Both of these facts can be verified with straightforward calculations.  Using these two facts, another straightforward calculation will show that $\nu'_{p}$ preserves the $\sigma$-eigenspace decomposition of $H^{1}(X_{p},\C)$. The second statement also follows from a direct computation.  \qed

\begin{notation} We denote by $\IAAI(\LXp)$ and $\IAAI(\Lambda)$ the sets of all involutive anti-isometries of $\LXp$ and $\Lambda$, respectively.
\end{notation}

\begin{definition}\label{FrameTypesIAAIclassesCorrespondence}\quad
We define the map $\pi_{0}(\FoR)\bigsqcup\pi_{0}(\Foa) \longrightarrow \P\IAAI(\Lambda)$
\begin{equation*}
\begin{array}{ccl}
\left(\,p,[i]\,\right) & \longmapsto & 
\left\{
\begin{array}{cll}
   \left[ \, i \circ \kappa^{*}_{p} \circ i^{-1} \, \right], & 
   \textnormal{if} \;\; p \in \PoR, & \textnormal{where $i$ is any frame over $p$}, \\
   \left[ \, i \circ \alpha^{*}_{p} \circ i^{-1} \, \right], & 
   \textnormal{if} \;\; p \in \Poa, &\textnormal{where $i$ is any frame over $p$}.
\end{array}
\right.
\end{array}
\end{equation*}
\end{definition}

\begin{definition}\label{TopTypesIAAIclassesCorrespondence}\quad  We also define
\begin{equation*}
\begin{array}{ccl}
\pi_{0}(\PoR)\bigsqcup\pi_{0}(\Poa) & \longrightarrow & \P\IAAI(\Lambda)/\P\Isom(\Lambda) \\
\left[\,p\,\right] & \longmapsto & 
\left\{
\begin{array}{cl}
   \left[ \, i \circ \kappa^{*}_{p} \circ i^{-1} \, \right], & \textnormal{if} \;\; p \in \PoR, \\
   \left[ \, i \circ \alpha^{*}_{p} \circ i^{-1} \, \right], & \textnormal{if} \;\; p \in \Poa.
\end{array}
\right.
\end{array}
\end{equation*}
\end{definition}

\begin{remark}\quad
The maps in Definitions \ref{FrameTypesIAAIclassesCorrespondence} and \ref {TopTypesIAAIclassesCorrespondence} are well-defined because $i \circ \kappa_{p}^{*}\circ i^{-1}$ and $i \circ \alpha_{p}^{*}\circ i^{-1}$ lie in the discrete subset $\IAAI(\Lambda)$ of $\IAAI(\Lambda\otimes_{\Zi}\C) \cong \IAAI(\C^{1,5})$, and hence remain constant as $p$ and $(\,p,[\,i\,]\,)$ vary within each connected component of $\PoR \bigsqcup \Poa$ and $\FoR\bigsqcup\Foa$ respectively.
\end{remark}


\subsection{Integral copies of $\RHfive$ in $\CHfive$}
\setcounter{theorem}{0}\setcounter{figure}{0}\setcounter{table}{0}

It can be readily checked that, for each $\chi\in\IAAI(\Lambda)$, the metric on $\Lambda$ restricts to a metric on the $\Z$-module $\Fix(\chi) \cong \Z^{6}$ of signature $(1+,5-)$.  Thus $\Fix(\chi)\otimes_{\Z}\Re$ $\cong$ $\Re^{1+,5-}$, and
\begin{equation*}
\begin{array}{ccc}
   \RH\left(\Fix(\chi)\otimes_{\Z}\Re\right) & \cong & \RHfive \\
   \cap                                      &       & \cap    \\
   \CH\left(\Lambda\otimes_{\Zi}\C\right)    & \cong & \CHfive
\end{array}
\end{equation*}
Hence, we may make the following

\begin{definition}\quad
A copy of $\RHfive \subset \CHfive$ is said to be \emph{integral} if it is of the form $\RH(\Fix(\chi)\otimes_{\Z}\Re)$ for some $\chi \in \IAAI(\Lambda)$.
\end{definition}


\subsection{``Real'' octics have ``real'' periods; ``antipodal'' octics have ``antipodal'' periods}
\setcounter{theorem}{0}\setcounter{figure}{0}\setcounter{table}{0}

Recall that, for any smooth $p \in \Po$,
\begin{equation*}
  \LXp \otimes_{\Zi}\C
  \cong \underset{\C^{1,5} = \C^{1+,5-}}{\underbrace{H^{1}_{\sigma=-\i}(X_{p},\C)}}
   =    \underset{(+)}{\underbrace{H^{1,0}_{\sigma=-\i}(X_{p},\C)}}
        \oplus \underset{(-----)}{\underbrace{H^{0,1}_{\sigma=-\i}(X_{p},\C)}}.
\end{equation*}

On the other hand, consider an ordered pair $(p,\nu_{p})$, where either $p\in\PoR$ and $\nu_{p} = \kappa_{p}$, or $p\in\Poa$ and $\nu_{p} = \alpha_{p}$.  Recall that $\nu'_{p}:H^{1}(X_{p},\C) \longrightarrow H^{1}(X_{p},\C)$ preserves both the Hodge decomposition and the $\sigma$-eigenspace decomposition.  Since $H^{1,0}_{\sigma=-\i}(X_{p},\C)$ is complex one-dimensional,  $H^{1,0}_{\sigma=-\i}(X_{p},\C) \in \CH\left(\LXp\otimes\C\right)$ is fixed by $[\nu_{p}] \in \P\IAAI(\LXp)$.  Hence, for a given framed smooth form $[\LXp\overset{i}{\rightarrow}\Lambda]$ over $p\in\PoR\bigsqcup\Poa$, and a fixed choice of $\nu_{p}$ ($=$ $\kappa_{p}$ or $\alpha_{p}$), the complex period $i(H^{1,0}_{\sigma=-\i}(X_{p},\C))$ $\in$ $\CHfive$ $=$ $\CH\left(\Lambda\otimes\C\right)$ is fixed by the projective class $[\chi_{\nu_{p}}]$ $=$ $[\,i \circ \nu^{*}_{p} \circ i^{-1}\,]$ $\in$ $\P\IAAI(\Lambda)$.  It now makes sense to introduce the following two definitions:

\begin{definition}\quad For each $[\,\chi\,] \in \P\IAAI(\Lambda)$, define $\RHfiveX$ to be the fixed point set of $[\,\chi\,]$ in $\CH(\Lambda\otimes_{\Zi}\C)\cong\CHfive$, i.e. $\RHfiveX$ $:=$  $\left\{\; [\,v\,]\in\CHfive \;\left\vert\;[\,\chi\,]([\,v\,]) = [\,v\,]\right. \;\right\}$.
\end{definition}

\begin{definition}  An element $x\in\CHfive$ is called a \emph{real period} if $x \in \RHfive_{[\chi_{\kappa_{p}}]}$, for some $p\in\PoR$.  An element $x\in\CHfive$ is called an \emph{antipodal period} if $x \in \RHfive_{[\chi_{\alpha_{p}}]}$, for some $p\in\Poa$.
\end{definition}

Let a representative $\chi \in [\,\chi\,] \in \P\IAAI(\Lambda)$ be fixed.  It is straightforward to see that we have the equality
\begin{equation*}
\RHfiveX \; = \;
\left\{\,[\,v\,]\in\CHfive\left\vert\;\exists\,v\in[\,v\,]\,\textnormal{\small with}\,\chi(v) = v \right.\,\right\}.
\end{equation*}
It is also easy to see that given any $[\,v\,] \in \RHfiveX$, the representative $v \in [\,v\,]$ that is fixed by the given $\chi$ is unique up to real scalar multiples.  This gives a canonical set-theoretic identification between $\RHfiveX$ and $\RH\left(\Fix(\chi)\,\otimes_{\Z}\,\Re\right) \cong \RHfive$.  The fixed point set $\RHfiveX$ is therefore canonically an integral copy of $\RHfive$ (hence its notation) and\, $\Stab_{\P\Isom\,\Lambda}\left(\RHfiveX\right)$ is isomorphic to a subgroup of $\Isom(\RHfive)$.  We see at once that the real and antipodal periods lie on integral copies of $\RHfive$ within $\CHfive$.


\subsection{The real period map and the Allcock-Carlson-Toledo construction of $\MoR$}
\setcounter{theorem}{0}\setcounter{figure}{0}\setcounter{table}{0}

The $G$-invariant complex period map $\p : \Fs \longrightarrow \CHfive$ was an important ingredient towards constructing the $\CHfive$ quotient structure for the moduli space $\Ms$ of stable complex binary octics.  We make use of it again to study the moduli space $\MsR$ of real binary octics.



\begin{definition}\quad
The \emph{real period map} is the map
\begin{equation*}
   \pR : \FoR \bigsqcup \Foa \longrightarrow \CHfive \times \P\IAAI(\Lambda)
\end{equation*}
defined by
\begin{equation*}
   \pR(\,p,[\,i\,]\,) \;\; := \;\;
   \left\{\begin{array}{cl}
      \left(\;\p(p,[i])\;,\;[\,i\circ\kappa_{p}^{*}\circ i^{-1}\,]\;\right), &
      \textnormal{if}\;\; (\,p\,,\,[\,i\,]\,) \in \FoR, \\
      \left(\;\p(p,[i])\;,\;[\,i\circ\alpha_{p}^{*}\circ i^{-1}\,]\;\right), &
      \textnormal{if}\;\; (\,p\,,\,[\,i\,]\,) \in \Foa. \\
   \end{array}\right.
\end{equation*}
\end{definition}

\begin{remark}\quad
The image of the real period map $\pR$ is naturally isomorphic (as real-analytic manifolds) to:
\begin{equation*}
   \Do \; := \; \bigsqcup_{[\chi]\in\P\IAAI(\Lambda)}\left(\RHfiveX-\H\right),
\end{equation*}
recalling that $\H\subset\CHfive$ is the collection of hyperplanes orthogonal to vectors in $\Lambda$ of squared norm $-2$.  Recall also that $\H$ is precisely the set of periods of singular octics.  Hereinafter, we regard $\Do$ as the codomain of $\pR$.
\end{remark}

\begin{definition}\quad
We let $\P\Gamma=\P\Isom(\Lambda)$ act on $\CHfive \times \P\IAAI(\Lambda)$ as follows: for $[\gamma]\in\P\Gamma$, and $\left(x,[\chi]\right) \in \CHfive \times \P\IAAI(\Lambda)$,
\begin{equation*}
     [\gamma]\cdot\left(\;x\;,\;[\chi]\;\right)
  := \left(\;\gamma(x)\;,\;[\,\gamma\circ\chi\circ\gamma^{-1}\,]\;\right).
\end{equation*}
This induces an action of $\P\Gamma$ on $\codomain(\pR)$ $=$ $\Do$ $=$ $\underset{[\chi]\in\P\IAAI(\Lambda)}{\bigsqcup}\left(\RHfiveX-\H\right)$.
\end{definition}

\begin{lemma}\label{RealPeriodGammaEquivariance}
The real period map is $\P\Gamma$-equivariant.
\end{lemma}

\begin{lemma}\label{RealPeriodGinvariance}
The real period map is $\GR$-invariant with respect to the action of $\GR$ on $\FoR$ and it is $\Ga$-invariant with respect to the action on $\Foa$.  In other words, it descends to a map, also denoted by $\pR$,
\begin{equation*}
  \pR \,:\, \left(\FoR/\GR\right) \bigsqcup \left(\Foa/\Ga\right)
         \, \longrightarrow \bigsqcup_{[\chi]\in\P\IAAI(\Lambda)}\RHfiveX.
\end{equation*}
Furthermore, the real period map $\pR$ restricts to a $\P\Gamma$-equivariant real-analytic diffeomorphism as follows:
\begin{equation*}
  \pR : \left(\FoR/\GR\right) \bigsqcup \left(\Foa/\Ga\right)
        \longrightarrow
        \Do := \bigsqcup_{[\chi]\in\P\IAAI(\Lambda)} \; \left(\RHfiveX - \H\right).
\end{equation*}
\end{lemma}



\begin{corollary}\label{IAAI:ClassesCount} \quad
The map defined in Definition \ref{TopTypesIAAIclassesCorrespondence}
\begin{equation*}
   \pi_{0}\left(\PoR\right)\,\bigsqcup\,\pi_{0}\left(\Poa\right)
   \longrightarrow
   \P\IAAI(\Lambda)/\P\Isom(\Lambda)
\end{equation*}
is surjective.  Consequently, the cardinality of $\P\IAAI(\Lambda)/\P\Isom(\Lambda)$ is at most seven.
\end{corollary}

The proofs of Lemmas \ref{RealPeriodGammaEquivariance}, \ref{RealPeriodGinvariance}, Corollary \ref{IAAI:ClassesCount}, as well as Proposition \ref{realperiodmapisisomorphism} involve unravelling the various definitions, the $G$-invariance and $\P\Gamma$-equivariance of the complex period map, and the fact that $\kappa_{p}^{*}$ commutes with elements of $\GR$ while $\alpha_{p}^{*}$ commutes with elements of $\Ga$.  Their complete proofs can be found in \cite{Chu:PhDthesis}, and they are straightforward adaptations of the proofs of the corresponding results in \cite{ACT:realcubics}.


\begin{lemma}\label{RealAntipIAAIClassesAreDisjoint}\quad
The images of $\pi_{0}(\PoR)$ and $\pi_{0}(\Poa)$ under the map in Definition \ref{TopTypesIAAIclassesCorrespondence} are disjoint in $\P\IAAI(\Lambda)/\P\Isom(\Lambda)$.
\end{lemma}
\proofoutline  This essentially follows from the observation that every octic form in $\PoR$ can deform to a nodal octic (i.e., a singular octic with one double root and no other singularities), whereas an octic in $\Poa$ can only deform to singular octics with at least two double points.  Recall that periods of nodal octics lie on the collection $\H \subset \CHfive$ of hyperplanes which are orthogonal complements of vectors in
$\Lambda$ of squared norm $-2$.  (See Section \ref{section:complexperiodmap}.)  By the preceding observation, we see that a copy of $\RHfive$ consisting of real periods must intersect $\H$ at smooth points of $\H$, whereas a copy of $\RHfive$ consisting of antipodal periods cannot intersect $\H$ at smooth points of $\H$.  \qed

\vskip 0.25cm
By Lemma \ref{RealAntipIAAIClassesAreDisjoint}, it makes sense to introduce the following:

\begin{definition} \quad
Let\; $\P\IAAI(\Lambda)^{\Re}/\P\Isom(\Lambda)$ and $\P\IAAI(\Lambda)^{\textnormal{antip}}/\P\Isom(\Lambda)$ be the images in $\P\IAAI(\Lambda)/\P\Isom(\Lambda)$ of $\pi_{0}(\PoR)$ and $\pi_{0}(\Poa)$, respectively, of the map
\begin{equation*}
  \pi_{0}\left(\PoR\right)\,\bigsqcup\,\pi_{0}\left(\Poa\right)
  \longrightarrow \P\IAAI(\Lambda)/\P\Isom(\Lambda),
\end{equation*}
as in Corollary \ref{IAAI:ClassesCount}.
\end{definition}

When we take into account the calculations in Section \ref{section:computations}, we furthermore have the following:

\begin{lemma}\label{RealIAAIClassesAreDisjoint}\quad
The values of $\PoRzero, \ldots, \PoRfour$ (\,considered as elements of $\pi_{0}(\PoR)$\,) under the map
$$\pi_{0}\left(\PoR\right)\,\bigsqcup\,\pi_{0}\left(\Poa\right)\longrightarrow\P\IAAI(\Lambda)/\P\Isom(\Lambda)$$
as in Definition \ref{TopTypesIAAIclassesCorrespondence} are pairwise distinct.
\end{lemma}
\proof The beginning of section \ref{section:computations} exhibits five involutive anti-isometries of $\Lambda$.  In section \ref{subsection:VinbergDiagrams}, it is shown that their fixed $\Z$-lattices have pairwise distinct Vinberg diagrams.  Hence, they represent five distinct $\P\Isom(\Lambda)$-conjugacy classes in $\P\IAAI(\Lambda)/\P\Isom(\Lambda)$.  Sections \ref{subsection:ChiZeroOneTwoThreeAreReal} and \ref{subsection:ChifourIsReal} show that all five are induced by real octics and identify their deformation types.
\qed

\begin{remark}\quad
We stress that Lemma \ref{RealIAAIClassesAreDisjoint} does not assert that $\PoRfourplus$ and $\PoRfourminus$ induce the same conjugacy class in $\P\IAAI(\Lambda)/\P\Isom(\Lambda)$; they may or may not.  However, this ambiguity does not pose a problem since our goal is just to describe the five connected components of $\MoR$ as abstract real hyperbolic quotients:  The complex linear change of variables $(x_{0},x_{1}) \longmapsto (\,\exp(\i\,\pi/8)\,x_{0},\,\exp(\i\,\pi/8)\,x_{1}\,)$ maps every $p(x_{0},x_{1}) \in \Po$ to $-p(x_{0},x_{1})$.  Consequently, even if the induced conjugacy classes in $\P\IAAI(\Lambda)/\P\Isom(\Lambda)$ of $\PoRfourplus$ and $\PoRfourminus$ are different, the respective real hyperbolic quotients will still be isomorphic.
\end{remark}

\begin{proposition}\label{realperiodmapisisomorphism}\quad
By further restricting the domain and codomain, and taking the quotient by $\P\Gamma$, the (\,$\P\Gamma$-equivariant\,) real period map
\begin{equation*}
  \pR : \left(\FoR/\GR\right) \bigsqcup \left(\Foa/\Ga\right)
        \longrightarrow
        \Do := \bigsqcup_{[\chi]\in\P\IAAI(\Lambda)} \; \left(\RHfiveX - \H\right)
\end{equation*}
descends to the following real-analytic orbifold isomorphism:
\begin{equation*}
  \MoR \bigsqcup \Moa
  \;\cong\;
  \P\Gamma \left\backslash \left(
  \textnormal{\small$\left(\FoRzero\sqcup\FoRone\sqcup\FoRtwo\sqcup\FoRthree\sqcup\FoRfourplus/\GR\right)$}
  \bigsqcup
  \textnormal{\small$\left(\Foa/\Ga\right)$}
  \right). \right.
\end{equation*}
In particular,
\begin{equation*}
  \MoR 
  \;\;\cong\;\;
  \P\Gamma \left\backslash \left(
  \FoRzero\sqcup\FoRone\sqcup\FoRtwo\sqcup\FoRthree\sqcup\FoRfourplus
  /\GR\right).\right.
\end{equation*}
\end{proposition}

\begin{figure}[b]
\begin{center}
\shortenedbigdiagram
\vskip 0.1cm
\caption{\small\it The decomposition of $\DoR$ and its projection to $\MoR$.}
\label{diagram:bigdiagram}
\end{center}
\end{figure}

The decomposition of $\DoR$ and its projection to $\MoR$ are illustrated in Figure \ref{diagram:bigdiagram}.  Combining Lemmas \ref{RealPeriodGinvariance}, \ref{RealAntipIAAIClassesAreDisjoint}, \ref{RealIAAIClassesAreDisjoint}, and Proposition \ref{realperiodmapisisomorphism}, we get the following

\begin{corollary}\label{ACT:MoR}\quad
Let $\chi_{0}, \ldots, \chi_{3}$, and  $\chi_{4}$ be any representatives of the conjugacy classes in $\P\IAAI(\Lambda)^{\Re}/\P\Isom(\Lambda)$ induced by $\PoRzero,\ldots,\PoRthree$, and $\PoRfourplus$, respectively.  Then, 
\begin{equation*}
   \MoRi \;\; \cong \;\; \P\GammaR_{i}\left\backslash \left(\RHfive_{\chi_{i}}-\H\right)\right.,
   \quad\;\;\textnormal{\it where}\quad
   \P\GammaR_{i}
   \;  = \;  \stab_{\P\Isom(\Lambda)}\left(\RHfive_{[\chi_{i}]}\right).
\end{equation*} 
Consequently,
\begin{equation*}
   \MoR \;\; = \;\; \bigsqcup_{i=0}^{4}\;\MoRi
   \;\; \cong \;\;
   \bigsqcup_{i=0}^{4}\;\P\GammaR_{i}\left\backslash \left(\RHfive_{\chi_{i}}-\H\right)\right..
\end{equation*} 
\end{corollary}

%
%
%

\subsection{The Allcock-Carlson-Toledo construction of $\MsR$}
\setcounter{theorem}{0}\setcounter{figure}{0}\setcounter{table}{0}

We will not give the full details of this construction.  Roughly speaking, it can be described as follows:
\begin{eqnarray*} 
                      \MsR
        &\;\;:=\;\;& \P(\PsR)\,/\;\PGLtwoR \,\; = \; \OsR \,/\, \PGLtwoR   \\
        & =       &  \left.\left(\OsRzero\cup\OsRone\cup\OsRtwo\cup\OsRthree\cup\OsRfour\right)\,\right/\,\PGLtwoR   \\
        &\cong  & 
                               \left.\left\{ \; 
                               \P\Gamma\,
                               \left\backslash\,\left(\FsRzero\cup\FsRone\cup\FsRtwo\cup\FsRthree\cup\FsRfourplus\right)\right.
                               \;\right\}\;\right/\,\GR \\
        &\cong  & 
                               \P\Gamma\,
                               \left\backslash\,\left\{\;
                               \left.\left(\FsRzero\cup\FsRone\cup\FsRtwo\cup\FsRthree\cup\FsRfourplus\right)\,\right/\GR
                               \;\right\}\right. \\
        & =     & \P\Gamma\,\backslash\,\Ks,
\end{eqnarray*}
where
\begin{equation*}
   \Ks := \longKs.
\end{equation*}
\noindent
Here, $\underset{0,1,2,3,4+}{\P\IAAI(\Lambda)}^{\Re}$ stands for the collection of all involutive anti-isometries of $\Lambda$ induced by smooth real binary octic forms of types $0$, $1$, $2$, $3$, and $4+$.  \; $\FsRzero$, $\FsRone$, $\FsRtwo$, $\FsRthree$, $\FsRfourplus$ are suitable completions of $\FoRzero$, $\FoRone$, $\FoRtwo$, $\FoRthree$, $\FoRfourplus$, respectively.  The quotient $\Ks$ is considered as a metric space, where the equivalence relation \;$\approx$\; on the disjoint union $\underset{[\chi]\in\underset{0,1,2,3,4+}{\P\IAAI(\Lambda)}^{\Re}}{\bigsqcup}\RHfive_{[\chi]}$ is defined so that we get a real-analytic homeomorphism
$\Ks \cong \left.\left(\FsRzero\cup\FsRone\cup\FsRtwo\cup\FsRthree\cup\FsRfourplus\right)\,\right/\GR$,
via standard arguments as in \cite{ACT:realcubics}.  

The Allcock-Carlson-Toledo construction is the metric space $\P\Gamma\,\backslash\,\Ks$, and it endows the moduli space $\MsR$ with a metric-space structure via the real-analytic homeomorphism mentioned above.  The crux of this construction therefore lies in explicitly describing the metric space $\Ks$ (i.e., how the disjoint copies $\RHfiveX$ are glued together) and the quotient $\P\Gamma\backslash\Ks$.

%% file: CompareStabilizers.tex
\section{Relationship between $\Stab_{\P\Isom\Lambda}\left(\RHfiveX\right)$ and $\P\Stab_{\Isom\Lambda}\left(\Fix\,\chi\right)$}
\setcounter{theorem}{0}\setcounter{figure}{0}\setcounter{table}{0}

In this section, we need to work simultaneously with projective equivalence classes of vectors, isometries and anti-isometries in various $\Z$-lattices and $\Zi$-lattices.  For the sake of clarity, we will use slightly more cumbersome notation such as $[\,v\,]_{\C}\in\CHfive$, $[\,A\,]_{\g} \in \P_{\g}\Isom\,\Lambda$ or $[\,A\,]_{\Z} \in \P_{\Z}\Isom(\Fix\,\chi)$ to indicate that the projectivization is done over $\C$, $\g = \Zi$ and $\Z$, respectively.


\subsection{Characterization of\, $\Stab_{\P_{\g}\Isom\,\Lambda}\left(\RHfive_{[\chi]}\right)$}
\setcounter{theorem}{0}\setcounter{figure}{0}\setcounter{table}{0}

Let $[\,\chi\,]_{\g} \in \P_{\g}\IAAI(\Lambda)$ be fixed.  Then,
\begin{equation*}
     \Stab_{\P_{\g}\Isom\,\Lambda}\left(\RHfive_{[\chi]}\right)
     \; := \;
     \left\{\; [\,A\,]_{\g}\in\P_{\g}\Isom\,\Lambda \;\left\vert\;
     [\,A\,]_{\g}\left(\RHfiveX\right)\subseteq\RHfiveX\right.\;\right\}.
\end{equation*}
Furthermore, let a representative $\chi \in [\,\chi\,]_{\g}$ be fixed.  Then, for $[\,A\,]_{\g}\in\P_{\g}\Isom(\Lambda)$, \vskip -0.2cm
\begin{equation*}
  \textnormal{\small$
  [\,A\,]_{\g} \in \Stab_{\P_{\g}\Isom\,\Lambda}\left(\RHfive_{[\chi]}\right)
  $}
  \;\Longleftrightarrow\;
  \textnormal{\small$
  \left\{\begin{array}{c}
       \textnormal{For each}\;A\in[\,A\,]_{\g},\;\textnormal{the following holds:}\\
       \textnormal{for each}\;
       [\,v\,]_{\C}\in\RHfiveX\;\textnormal{and}\;
       v\in[\,v\,]_{\C}\;\textnormal{with}\;\chi(v)=v, \\
       \exists\;\textnormal{unique}\;\beta\in\C^{*}\;\textnormal{with}\;|\beta|=1\;\textnormal{and}\;
       \chi(A(v)) = \beta A(v).
   \end{array}\right.
   $}
\end{equation*}\vskip -0.2cm

\begin{remark}\quad
The uniqueness (once the representatives $A \in [\,A\,]_{\g}$ and $\chi \in [\,\chi\,]_{\g}$ are fixed) and unimodularity of $\beta$ above are clear.  Since both $A$ and $\chi$ preserve primitiveness of lattice vectors, we see that $\beta$ is in fact a unit Gaussian integer whenever $v \in \Fix(\chi)$ is primitive in $\Lambda$.  If $v \in \Fix(\chi)$ is only primitive in the $\Z$-lattice $\Fix(\chi)$, but not in $\Lambda$, then $v = (1+\i)w$, for some $w$ primitive in $\Lambda$.  It can be readily shown that $\chi(A(v)) = \beta\,A(v)$ implies $\chi(A(w)) = \i\,\beta\,A(w)$.  $\Lambda$-primitiveness of $w$ then again shows that $\beta$ must be a unit Gaussian integer.
\end{remark}

\begin{lemma}\quad Let $\chi \in \IAAI(\Lambda)$ be given.  Let $A \in \Isom(\Lambda)$ be such that $[\,A\,]_{\g} \in \Stab_{\P_{\g}\Isom\,\Lambda}\left(\RHfive_{[\chi]}\right)$.   Then, there exists a unique $\beta \in \C^{*}$ such that $\chi(A(v)) = \beta A(v)$, for all $v \in \Fix(\chi)\otimes_{\Z}\Re$.  Furthermore, $\beta$ is in fact a unit Gaussian integer.
\end{lemma}
\proofoutline   From the preceding remark, we know that for each given $v \in \Fix(\chi)\otimes_{\Z}\Re$, there exists a unique unimodular $\beta \in \C^{*}$ such that $\chi(A(v)) = \beta\,A(v)$.  Furthermore, $\beta$ is a unit Gaussian integer whenever $v$ is primitive in the $\Z$-lattice $\Fix(\chi)$.  So, it remains to show only that $\beta$ is in fact the same for all $v \in \Fix(\chi)\otimes_{\Z}\Re$.   For this, let $b_{1},\ldots,b_{6}$ be a $\Z$-basis for $\Fix(\chi)$ and let $v \in \Fix(\chi)\otimes_{\Z}\Re$.  Then, there exist unit Gaussian integers $\beta_{1},\ldots,\beta_{6}\in\Zi$ such that $\chi(A(b_{k})) = \beta_{k}A(b_{k})$, unique for each $k = 1, \ldots, 6$.  Also, there exists unique unimodular $\beta \in \C^{*}$ such that $\chi(A(v)) = \beta\,A(v)$.  Now, recall that $\Fix(\chi)\otimes_{\Z}\Re$ is a totally real subspace of $\C^{1,5}$.   In particular, $b_{1},\ldots,b_{6}$ are linearly independent over $\C$.  This observation and a simple calculation show that we must in fact have $\beta_{1} = \cdots =\beta_{6} =\beta$, which completes the proof.  \qed


\begin{proposition}\label{typeIandII}\quad
Let $[\,\chi\,]_{\g} \in \P_{\g}\IAAI(\Lambda)$ be fixed and fix also a representative $\chi \in [\,\chi\,]_{\g}$.  Then, for each $[\,A\,]_{\g} \in \Stab_{\P_{\g}\Isom\Lambda}\left(\RHfiveX\right)$, exactly one of the following holds:
\begin{eqnarray*}
   \textnormal{either} & \;\; &
   \exists \; A \in [\,A\,]_{\g},\,\textnormal{unique up to sign, such that}\;
   A(\Fix(\chi)) \subseteq \Fix(\chi), \\
   \textnormal{or} & \;\; &
   \exists \; A \in [\,A\,]_{\g},\,\textnormal{unique up to sign, such that}\;
   A(\Fix(\chi)) \subseteq \Fix(\i\,\chi).
\end{eqnarray*}
\end{proposition}
\proof We already know that, for an arbitrary representative $A \in [\,A\,]_{\g}$, we have $A(\Fix\,\chi) = \Fix(\beta\,\chi)$, where $\beta$ is one of the four unit Gaussian integers.  The Proposition thus trivially follows from observing what the other associates of $A$ are doing to $\Fix(\chi)$, as shown in Table \ref{table:AFixChi}.  \qed

\begin{table}
\begin{center}
\begin{tabular}{|c||c|c|c|}
\hline
$A(\Fix\,\chi)$ & $-A(\Fix\,\chi)$ & $\i\,A(\Fix\,\chi)$ & $-\i\,A(\Fix\,\chi)$ \\ \hline\hline
$\Fix(\chi)$     & $\Fix(\chi)$      & $\Fix(-\chi)$        & $\Fix(-\chi)$         \\ \hline 
$\Fix(-\chi)$    & $\Fix(-\chi)$     & $\Fix(\chi)$         & $\Fix(\chi)$          \\ \hline 
$\Fix(\i\,\chi)$  & $\Fix(\i\,\chi)$  & $\Fix(-\i\,\chi)$    & $\Fix(-\i\,\chi)$     \\ \hline 
$\Fix(-\i\,\chi)$ & $\Fix(-\i\,\chi)$ & $\Fix(\i\,\chi)$     & $\Fix(\i\,\chi)$      \\
\hline 
\end{tabular}
\vskip 0.1cm
\caption{}
\label{table:AFixChi}
\end{center}
\end{table}

\begin{remark}\quad
Note that if $A(\Fix\,\chi) = \Fix(\chi)$, then $A$ also preserves each of $\Fix(\i\,\chi)$, $\Fix(-\chi)$, $\Fix(-\i\,\chi)$.  On the other hand, if $A(\Fix\,\chi) = \Fix(\i\,\chi)$, then $A$ maps $\Fix(\i\,\chi)$ to $\Fix(-\chi)$, $\Fix(-\chi)$ to $\Fix(-\i\,\chi)$, and $\Fix(-\i\,\chi)$ to $\Fix(\chi)$. Hence, we have the following:
\end{remark}

\begin{proposition}\label{CharacterizingStab}\quad
Let $\chi \in [\,\chi\,]_{\g}$ be fixed.  Then, the stabilizer $\Stab_{\P_{\g}\Isom\,\Lambda}\left(\RHfive_{[\chi]}\right)$ can be characterized as follows:
\begin{equation*}
        \Stab_{\P_{\g}\Isom\,\Lambda}\left(\RHfive_{[\chi]}\right)
   =  \textnormal{\large$\P_{\Z}$}
        \Stab_{\Isom\,\Lambda}
        \left(\!\!\begin{array}{c}
                      \\ \Fix(\chi)\cup\Fix(\i\chi)\cup\Fix(-\chi)\cup\Fix(-\i\chi) \\ \\
        \end{array}\!\!\right).
\end{equation*}
\end{proposition}

We seek an even more algebraically transparent expression for $\Stab_{\P_{\g}\Isom\,\Lambda}(\RHfiveX)$ in terms of\, $\P_{\Z}\Stab_{\Isom\,\Lambda}\left(\Fix\,\chi\right)$.

\begin{definition}\quad
$[\,A\,]_{\g} \in \Stab_{\P_{\g}\Isom\Lambda}\left(\RHfiveX\right)$ is said to be of \emph{type I} if there exists $A \in [\,A\,]_{\g} \in \P_{\g}\Isom\Lambda$ such that $A(\Fix\chi)\subset\Fix(\chi)$, and it is said to be of \emph{type II} if there exists $A \in [\,A\,]_{\g} \in \P_{\g}\Isom\Lambda$ such that $A(\Fix\chi)\subset\Fix(\i\chi)$.
\end{definition}

\begin{remark}\quad
Note that $[\,A\,]_{\g} \in \Stab_{\P_{\g}\Isom\Lambda}\left(\RHfiveX\right)$ is either of type I or type II by Proposition \ref{typeIandII}.  Some simple calculations will furthermore show:
\end{remark}

\begin{lemma}\label{lemma:CompositionAndInverse}\quad
\begin{enumerate}
\item  If two elements in $\Stab_{\P_{\g}\Isom\Lambda}\left(\RHfiveX\right)$ are of the same type, then their composition is an element of type I.
\item  If two elements in $\Stab_{\P_{\g}\Isom\Lambda}\left(\RHfiveX\right)$ are of different types, then their composition in either order is of type II.  
\item  Taking inverses in $\Stab_{\P_{\g}\Isom\Lambda}\left(\RHfiveX\right)$ preserves types.
\end{enumerate}
\end{lemma}
It is already clear that either
\begin{equation*}
\Stab_{\P_{\g}\Isom\Lambda}\left(\RHfiveX\right) = \P_{\Z}\Stab_{\Isom\Lambda}(\Fix\,\chi),
\quad\quad\textnormal{or}\quad\quad
\dfrac{\Stab_{\P_{\g}\Isom\Lambda}\left(\RHfiveX\right)}{\P_{\Z}\Stab_{\Isom\Lambda}(\Fix\,\chi)} \cong \Z/2\Z.
\end{equation*}
We will next show, under the further assumption that $\P_{\Z}\Stab_{\Isom\,\Lambda}(\Fix\,\chi)$ is a reflection group, that the following short exact sequence
\begin{equation*}
   1 \,\longrightarrow\, \P_{\Z}\Stab_{\Isom\,\Lambda}(\Fix\,\chi) 
      \,\longrightarrow\, \Stab_{\P_{\g}\Isom\,\Lambda}\left(\RHfiveX\right)
      \,\longrightarrow\, \Z/2\Z
      \,\longrightarrow\, 1
\end{equation*}
is split (Proposition \ref{Stabilizer:semidirectproduct}).
We start with the following general result, the proof of which
follows from standard arguments and we therefore omit:

\begin{proposition}\label{general:semidirectproduct}\quad
Let $G$ be a discrete subgroup of $\Isom(\RH^{n})$.  Suppose that $H$ is a normal subgroup of $G$ which is generated by reflections.  ($H$ need not be the full reflection subgroup of $G$.)  Fix a fundamental domain $P$ of $H$, and let $K := \{\;g \in G \,\;|\;\, g \cdot P = P \;\}$.  Then, $G \; = \; H \rtimes K$, where the action of $K$ on $H$ is, as usual, by conjugation.
\end{proposition}





\begin{remark}\quad
$K$ in Proposition \ref{general:semidirectproduct} is a subgroup of the symmetry group of the fundamental domain $P$.  $K$ may be trivial even if the symmetry group of $P$ is not.  Obviously, if $K$ is trivial, then $G = H$.
\end{remark}

Recall that $\P_{\Z}\Stab_{\Isom\,\Lambda}\left(\Fix\,\chi\right)$ is a normal subgroup of $\Stab_{\P_{\g}\Isom\,\Lambda}\left(\RHfiveX\right)$ of index two or one, depending on whether or not there are elements of type II.  Using Proposition \ref{general:semidirectproduct}, we now obtain the following expression for $\Stab_{\P_{\g}\Isom\,\Lambda}\left(\RHfiveX\right)$ in terms of $\P_{\Z}\Stab_{\Isom\,\Lambda}(\Fix\,\chi)$:

\begin{proposition}\label{Stabilizer:semidirectproduct}\quad
Suppose $\P_{\Z}\Stab_{\Isom\Lambda}(\Fix\,\chi)$ is generated by reflections.  Then, exactly one of the following holds:
\begin{itemize}
     \item  $\Stab_{\P_{\g}\Isom\Lambda}\left(\RHfiveX\right)$ has no elements of type II, in which case,
               \begin{equation*}
                    \Stab_{\P_{\g}\Isom\Lambda}\left(\RHfiveX\right) = \P_{\Z}\Stab_{\Isom\Lambda}(\Fix\,\chi),
               \end{equation*}
     \item  $\Stab_{\P_{\g}\Isom\Lambda}\left(\RHfiveX\right)$ contains elements of type II, in which case, 
               the fundamental domain of the group action
               $\P_{\Z}\Stab_{\Isom\,\Lambda}\left(\Fix\,\chi\right)\acts\RHfive$ admits a $(\Z/2\Z)$-symmetry, and
               via its norm-preserving action on the roots of $\P_{\Z}\Stab_{\Isom\,\Lambda}\left(\Fix\,\chi\right)$,
               this $(\Z/2\Z)$-symmetry induces an order-two element
               $[\,T\,] \in \Stab_{\P_{\g}\Isom\,\Lambda}\left(\RHfiveX\right)$
               of type II such that
               \begin{eqnarray*}
                                  \Stab_{\P_{\g}\Isom\Lambda}\left(\RHfiveX\right)
                    &   =    & \P_{\Z}\Stab_{\Isom\Lambda}(\Fix\,\chi) \rtimes \langle\,[\,T\,]\,\rangle \\
                    &\cong& \P_{\Z}\Stab_{\Isom\Lambda}(\Fix\,\chi) \rtimes \left(\Z/2\,\Z\right).
               \end{eqnarray*}
\end{itemize}
\end{proposition}

\begin{remark}\quad
Any representative $T \in \Isom(\Lambda)$ of the type II and order-two element $[\,T\,] \in \Stab_{\P_{\g}\Isom\,\Lambda}\left(\RHfiveX\right)$ maps $\Fix(\chi)$ to $\Fix(\i\,\chi)$, rather than back to $\Fix(\chi)$ itself.  $T$ induces an action on $\RHfiveX \cong \RH(\Fix(\chi)\otimes_{\Z}\Re)$ by identifying $\RH(\Fix(\i\,\chi)\otimes_{\Z}\Re)$ with $\RH(\Fix(\chi)\otimes_{\Z}\Re)$ via scalar multiplication by $(1-\i)/\sqrt{2}$; more explicitly,
\begin{equation*}
\begin{array}{rcl}
\RH(\Fix(\i\,\chi)\otimes_{\Z}\Re) & \longrightarrow & \RH(\Fix(\chi)\otimes_{\Z}\Re)       \\
  \left[\,w\,\right]                            & \longmapsto     & \left[\,\frac{1-\i}{\sqrt{2}}\,w\,\right]
\end{array}.
\end{equation*}
This identification is canonical due to the following observation:
\begin{equation*}
                               \i\chi(w) = w
   \;\;\Longleftrightarrow\;\; \left(\frac{1+\i}{\sqrt{2}}\right)\left(\frac{1+\i}{\sqrt{2}}\right)\chi(w) = w
   \;\;\Longleftrightarrow\;\; \chi\left(\frac{1-\i}{\sqrt{2}}w\right) = \left(\frac{1-\i}{\sqrt{2}}\right)w.
\end{equation*}
We emphasize that while $T$ preserves the $\Re$-span of $\Fix(\chi)$ via the above canonical induced action, it fails to preserve the $\Z$-lattice $\Fix(\chi)$ itself due to the occurrence of the $1/\sqrt{2}$ factor above.
\end{remark}



\subsection{A sufficient condition for the nonexistence of isometries of type II}
\setcounter{theorem}{0}\setcounter{figure}{0}\setcounter{table}{0}

Note that $\P_{\Z}\Stab_{\Isom\,\Lambda}\left(\Fix\,\chi\right)$ is merely the subgroup of the isometry group $\P_{\Z}\Isom(\Fix\,\chi)$ of the abstract $\Z$-lattice $\Fix(\chi)$ consisting of elements that extend to an action on the whole $\Zi$-lattice $\Lambda$.  In the case where $\P_{\Z}\Stab_{\Isom\Lambda}(\Fix\,\chi)$ is a reflection group and $\Stab_{\P_{\g}\Isom\,\Lambda}\left(\RHfiveX\right)$ contains type II elements, we see that we have the following commutative diagram:
\begin{diagram}
\P_{\Z}\Stab_{\Isom\,\Lambda}\left(\Fix\,\chi\right)&\rInto^{\mbox{}\quad\quad\mbox{}}&
\P_{\Z}\Stab_{\Isom\,\Lambda}\left(\Fix\,\chi\right) \rtimes \langle[\,T\,]\rangle
& \mbox{} \;\, = \;\, \mbox{} & \Stab_{\P_{\g}\Isom\,\Lambda}\left(\RHfiveX\right) \\
\dInto & & \dInto & & \\
\P_{\Z}\Isom\left(\Fix\,\chi\right) & \rInto & \P_{\Z}\Isom\left(\Fix\,\chi\right)\rtimes \langle[\,T\,]\rangle &&
\end{diagram}
where 
$[\,T\,] \in \Stab_{\P_{\g}\Isom\,\Lambda}\left(\RHfiveX\right)$ is an element of type II and order two.  Proposition \ref{general:semidirectproduct} therefore implies the following:

\begin{corollary}\label{NecessaryCondn1:ExistenceTypeII}\quad
Suppose $\P_{\Z}\Isom(\Fix\,\chi)$ is generated by reflections, and suppose one of the following conditions holds:
\begin{itemize}
     \item   The fundamental domain of $\P_{\Z}\Isom(\Fix\,\chi)$ admits no $(\Z/2\Z)$-symmetries.
     \item   It admits $(\Z/2\Z)$-symmetries, but none of them induces an order-two element of\;
                $\Stab_{\P_{\g}\Isom\Lambda}\left(\RHfiveX\right)$ of type II.
\end{itemize} \vskip -0.2cm
Then, $\Stab_{\P_{\g}\Isom\Lambda}\left(\RHfiveX\right)$ in fact has no elements of type II, and 
\begin{equation*}
     \Stab_{\P_{\g}\Isom\Lambda}\left(\RHfiveX\right) = \P_{\Z}\Stab_{\Isom\Lambda}(\Fix\,\chi).
\end{equation*}
\end{corollary}


\begin{remark}\quad
In \cite{Chu:PhDthesis}, $\Stab_{\P_{\g}\Isom\Lambda}\left(\RHfiveX\right)$ was mistakenly identified with $\P_{\Z}\Stab_{\Isom\Lambda}(\Fix\,\chi)$ in general, which need not be the case as we saw in this section.   The main results stated there are nonetheless correct, since for the specific cases therein (i.e. $\chi = \chizero, \chione, \chitwo, \chifour$), the above equality indeed holds.  
\end{remark}


%% file: InducedPermutation.tex
\section{Distinguishing the Deformation Types}
\label{chapter:IdentifyingTopTyesOfFixedLattices}
\setcounter{theorem}{0}\setcounter{figure}{0}\setcounter{table}{0}


In this section, we describe a strategy to identify the deformation types of the real octics that give rise to the involutive anti-isometries of $\Lambda$.


\subsection{The isomorphism $\O\!\left(\textnormal{\small$\Lambda/(1+\i)\,\Lambda$},q\right)$ $\cong$ $S_{8}$}
\setcounter{theorem}{0}\setcounter{figure}{0}\setcounter{table}{0}


Let $h$ be the $\Zi$-valued inner product of $\Lambda$ and $Q$ be the associated $\Zi$-valued quadratic form.  $Q$ is ``even-valued,'' and $\frac{1}{2}Q$ is thus a well-defined $\Z$-valued function on $\Lambda$.  On the other hand, $\Zi / (1+\i)\,\Zi \cong \F_{2}$, as rings (hence as fields), where $\F_{2}$ denotes the field with two elements.  $V := \Lambda/(1+\i)\Lambda$ is a six-dimensional $\F_{2}$-vector space.  The $\F_{2}$-valued function $q$ on $V$ defined by $x \overset{q}{\longmapsto} \frac{1}{2}Q(x)\!\mod (1+\i)$ is an $\F_{2}$-valued quadratic form on $V$.  It turns out that the orthogonal group $\O(V,q)$ is isomorphic to $S_{8}$, the symmetric group on eight objects.

We will not give complete proofs of the above assertions but refer the reader to \cite{Chu:PhDthesis} and \cite{AtlasFG}.  However, we give an intuitive description of the isomorphism between $\O(V,q)$ and $S_{8}$.

Since $\dim_{\F_{2}}(V) = 6$, we immediately see that the cardinality of $V$ is $2^{6}=64$.  Let $P_{8} := \{1,\ldots,8\}$.  Here, we simply consider $P_{8}$ to be any set of the eight distinct objects.  It turns out that, as a set, $V$ is in one-to-one correspondence with
\begin{equation*}
   W := \left.\left\{\begin{array}{c}
   \textnormal{even-cardinality} \\ \textnormal{subsets of $P_{8}$}
   \end{array}\right\}\right/ 
   \left\{B \;\, \sim
      \begin{array}{c}\textnormal{complement}\\
      \textnormal{of}\;\,B\;\,\textnormal{in $P_{8}$}\end{array}
   \right\}.
\end{equation*}
Each element of $W$ can be considered as a pair of even-cardinality subsets of $P_{8}$, where the two subsets in each such pair are complements of each other.  The cardinality of $W$ is also $64$.  The $\F_{2}$-valued quadratic form on $V$ corresponds to the $\F_{2}$-valued function on $W$ given by: 
\begin{equation*}
\begin{array}{rcl}
     W & \longrightarrow & \F_{2} \\
      s & \longmapsto      & \frac{1}{2}\left(\textnormal{cardinality of}\;\,s\right)
\end{array}.
\end{equation*}
Furthermore, elements of $\O(V,q)$ correspond to maps from $P_{8}$ to itself which preserve the cardinality of every even-cardinality subset of $P_{8}$.  Such a map is just a permutation of $P_{8}$, namely, an element of the symmetric group $S_{8}$.  It turns out that this map $\O(V,q) \longrightarrow S_{8}$ is an isomorphism of groups.  We denote its inverse by $\Phi:S_{8}\longrightarrow\O(V,q)$.


\subsection{Two useful invariants of involutions in $\Seight$}
\setcounter{theorem}{0}\setcounter{figure}{0}\setcounter{table}{0}

Recall that the eight distinct roots of a smooth real binary octic are preserved as a set by complex conjugation $\conjugationCPone$ on $\CPone$.  The collection $P_{8}$ of roots comprises a number $2n \in \{0,2,4,6,8\}$ of real points (lying on $\RPone = \Re \cup \{\infty\} \subseteq \CPone$) together with a number $(8 - 2n)/2$ of complex conjugate pairs.  The number $2n$ determines the deformation type of a real binary octic.

On the other hand, note that when $\conjugationCPone$ is restricted to the collection $P_{8}$ of the eight distinct roots of a real binary octic, it becomes an order-two permutation on $P_{8}$.  Table \ref{table:OcticTypesVsCycleStructures} shows the one-to-one correspondence between the deformation types of octics with the cycle structures of $\conjugationCPone|_{P_{8}}$.
\begin{table}
\begin{center}
\begin{tabular}{|c|c|c|}
   \hline
   $\begin{array}{c}\textnormal{\small Type} \\ \textnormal{\small of octic}\end{array}$
    & $2n$ & cycle structure of $\conjugationCPone|_{P_{8}}$ \\ \hline\hline
    0   & 8    & $(1)(2)(3)(4)(5)(6)(7)(8)$ \\ \hline
    1   & 6    & $(1)(2)(3)(4)(5)(6)(78)$   \\ \hline
    2   & 4    & $(1)(2)(3)(4)(56)(78)$     \\ \hline
    3   & 2    & $(1)(2)(34)(56)(78)$       \\ \hline
    4   & 0    & $(12)(34)(56)(78)$         \\ \hline
\end{tabular}
\vskip 0.1cm
\caption{}
\label{table:OcticTypesVsCycleStructures}
\end{center}
\end{table}
\vskip 0.2cm

Of course the cycle structure of $\conjugationCPone|_{P_{8}}$ determines an element (more precisely, a conjugacy class) in $\Seight$.  Now we make the following observations:
\begin{itemize}
   \item  $\conjugationCPone$ induces an involutive
             antiholomorphic diffeomorphism
             $\kappa: X_{p} \longrightarrow X_{p}$ on the 4-sheeted cyclic
             cover $X_{p} \longrightarrow \CPone$ branched over the roots
             of a smooth real binary octic form $p$.  $\kappa$ in turn induces
             an involutive anti-isometry on the
             $\Zi$-lattice $\LXp \cong \Lambda$.
   \item  $(V,q) := (\,\textnormal{\small$\Lambda\,/\,(1+\i)\Lambda$}\,,\,q\,)$
             is an orthogonal space over $\Zi/(1+\i) \cong \F_{2}$ such that
             $\O(V,q) \cong S_{8}$.
   \item  The above ``abstract'' isomorphism $\O(V,q) \cong \Seight$ is
             geometrically realized by the permutation of the eight ramification
             points of the branched cover of $X \longrightarrow \CPone$.
             This fact is an immediate consequence of the fact that the
             monodromy group $\P\Gamma$ $=$ $\P\Isom(\Lambda)$ is
             generated by transposing pairs of roots by ``continuous half
             turns.''  See \cite{MY:complex8points}.             
   \item  An involutive anti-isometry of $\Lambda$ descends to an
             involutive isometry of $(V,q)$ (because complex conjugation
             on $\Zi$ descends to the identity on $\Zi/(1+\i)\Zi \cong \F_{2}$).
\end{itemize}
The above observations show the following: Given $\chi \in \IAAI(\Lambda)$, we can determine the
deformation type of the real binary octic that gives rise to $\chi$ in the first place by determining the element (or conjugacy class) in $\Seight \cong \O(V,q)$ that $\chi$ descends to.  In order to do this, it is sufficient to examine two invariants:

\begin{lemma} \label{LemmaTwoInvariants}\quad
Let $\Phi : \Seight \longrightarrow \O(V,q)$ be the isomorphism (unique up to conjugacy) constructed earlier.
Then, the invariants $\dim_{\F_{2}}\Fix(\Phi(\tau_{i}))$ and the number of norm-one vectors in $\Fix(\Phi(\tau_{i}))$ of the various cycle structures are as shown in Table \ref{table:CycleStructureWithInvariants}.
\end{lemma}

\begin{table}
\begin{center}
\begin{tabular}{|c|c|c|c|}
   \hline
   Type & cycle structure of $\conjugationCPone|_{P_{8}}$ &
   $\dim_{\F_{2}}\Fix(\Phi(\tau_{i}))$             &
   $\begin{array}{c}
     \textit{\small number of norm-one} \\
     \textit{\small vectors in $\Fix(\Phi(\tau_{i}))$}
    \end{array}$ \\
   \hline\hline
   0 & $\tau_{8} = (1)(2)(3)(4)(5)(6)(7)(8)$ & 6 & 28 \\ \hline
   1 & $\tau_{6} = (1)(2)(3)(4)(5)(6)(78)$   & 5 & 16 \\ \hline
   2 & $\tau_{4} = (1)(2)(3)(4)(56)(78)$     & 4 &  8 \\ \hline
   3 & $\tau_{2} = (1)(2)(34)(56)(78)$       & 3 &  4 \\ \hline
   4 & $\tau_{0} = (12)(34)(56)(78)$         & 4 &  4 \\ \hline
\end{tabular}
\vskip 0.1cm
\caption{}
\label{table:CycleStructureWithInvariants}
\end{center}
\end{table}

\proofoutline  Let $P_{8} = \{1,2,\ldots,8\}$.  Recall that norm-one vectors in $V$ correspond to cardinality-two subsets of $P_{8}$.  The computations for all the cases are similar; we show only those for $\tau_{6}$:  The number of even-cardinality subsets of $P_{8}$ fixed by $\tau_{6} = (1)(2)(3)(4)(5)(6)(78)$ is given by
{\small\begin{equation*}
   2 \times \left(\;   \left(\begin{array}{c}6\\0\end{array}\right)
                     + \left(\begin{array}{c}6\\2\end{array}\right)
                     + \left(\begin{array}{c}6\\4\end{array}\right)
                     + \left(\begin{array}{c}6\\6\end{array}\right)
            \;\right)
   \; = \;
   2 \times \left(\; 1 + 15 + 15 + 1 \;\right) \; = \; 2 \times 32.  
\end{equation*}}
Hence, $\dim_{\F_{2}}\Fix(\Phi(\tau_{6})) = \log_{2}\left(\dfrac{2 \times 32}{2}\right) = \log_{2}(2^{5}) = 5$.  And,
\begin{eqnarray*}
  \textnormal{\small$
  \left\{\begin{array}{c}
     \textnormal{\small the number of} \\
     \textnormal{\small norm-one} \\
     \textnormal{\small vectors} \\
     \textnormal{\small in $\Fix(\Phi(\tau_{6}))$}
  \end{array}\right\}$}
  & = &
  \textnormal{\small$
  \left\{\begin{array}{c}
     \textnormal{\small the number of} \\
     \textnormal{\small cardinality-two} \\
     \textnormal{\small subsets} \\
     \textnormal{\small preserved by $\tau_{6}$}
  \end{array}\right\}$}
  \; = \;
  1 + 
  \textnormal{\small$
  \left\{\begin{array}{c}
     \textnormal{\small the number of all} \\
     \textnormal{\small cardinality-two} \\
     \textnormal{\small subsets of $\{1,\ldots,6\}$}
  \end{array}\right\}$} \\
  & = &
  1 + \left(\begin{array}{c}6 \\ 2\end{array}\right)
  \; = \;
  16.
\end{eqnarray*} \vskip -0.5cm \qed

\begin{remark}\quad
The antipodal map $\CPone\overset{\alpha}{\longrightarrow}\CPone$ permutes the roots of a smooth antipodal octic in the same way as complex conjugation $\CPone\overset{\kappa}{\longrightarrow}\CPone$
does the roots of a smooth real octic of type $4$ (i.e. the roots are four complex conjugate pairs).   The cycle structure for both is $(12)(34)(56)(78)$.  Hence, Lemma \ref{LemmaTwoInvariants} is insufficient to distinguish an antipodal octic from a real octic of type $4$.  To achieve this, we will need the idea of the proof of Lemma \ref{RealAntipIAAIClassesAreDisjoint} instead.  See section \ref{subsection:ChifourIsReal}.
\end{remark}


%% file: ComputationalResults.tex
\section{Computations}\label{section:computations}
\setcounter{theorem}{0}\setcounter{figure}{0}\setcounter{table}{0}


Define the map $\chitwo:\Lambda\longrightarrow\Lambda$ by
\begin{equation*}
     \chitwo\textnormal{\scriptsize$\zvector$} \; := \; \textnormal{\scriptsize$\chitwomatrix$}.
\end{equation*}
Clearly, $\chitwo \in \IAAI(\Lambda)$.  Next, define
\begin{equation*}
\begin{array}{lll}
     \Azero  :=  \textnormal{\scriptsize$\Azeromatrix$}, &
     \Aone   :=  \textnormal{\scriptsize$\Aonematrix$},  &
     \Atwo    :=  \textnormal{\scriptsize$\Atwomatrix$},  
\end{array}
\end{equation*}
\begin{equation*}
\begin{array}{cccccccccc}
     \Athree &:=& \textnormal{\scriptsize$\Athreematrix$},
     & \mbox{}\quad\mbox{} &
     \Afour   &:=& \textnormal{\scriptsize$\Afourmatrix$}.
\end{array}
\end{equation*}
Consider $\Azero, \ldots, \Afour$, as $\Zi$-linear endomorphisms on $\Lambda$, via $v \mapsto A_{i}\cdot v$.  Define
\begin{equation*}
     \chizero := \Azero  \circ \chitwo, \;\;\;
     \chione  := \Aone   \circ \chitwo, \;\;\;
     \chitwo   := \Atwo    \circ \chitwo \, = \, \chitwo, \;\;\;
     \chithree := \Athree \circ \chitwo, \;\;\;
     \chifour  := \Afour   \circ \chitwo.
\end{equation*}
It is straightforward to verify that $\Azero, \ldots, \Afour \in \Isom(\Lambda)$, and $\chizero, \ldots, \chifour \in \IAAI(\Lambda)$.  As the notation suggests, $\chizero, \ldots, \chifour$ correspond to real binary octics of types $0, \ldots, 4$ respectively.  Appealing to the theory developed in the preceding sections, we now present a series of straightforward computational results in the following subsections which will establish this correspondence.  We will omit the details of these computations due to their routine nature.


\subsection{$\Z$-bases for the fixed lattices of $\chizero,\ldots,\chifour$}
\label{subsection:FixedBases}
\setcounter{theorem}{0}\setcounter{figure}{0}\setcounter{table}{0}

The column vectors of the following matrices form respectively $\Z$-bases for the
fixed $\Z$-lattices of the anti-involutions $\chizero, \ldots, \chifour$: 
               \begin{equation*}
               \begin{array}{ll}
                    \basiszero  :=\textnormal{\tiny$\Basiszeromatrix$},  &
                    \basisone   :=\textnormal{\tiny$\Basisonematrix$},   \\ \\
                    \basistwo    :=\textnormal{\tiny$\Basistwomatrix$},   &
                    \basisthree :=\textnormal{\tiny$\Basisthreematrix$}, \\ \\
                    \basisfour   :=\textnormal{\tiny$\Basisfourmatrix$}.  &
               \end{array}
               \end{equation*}
               

\subsection{The induced integral quadratic forms on $\Fix(\chizero),\ldots,\Fix(\chifour)$}
\label{subsection:FixedLattices}
\setcounter{theorem}{0}\setcounter{figure}{0}\setcounter{table}{0}

These are determined by inner product matrices of $\basiszero, \ldots, \basisfour$, which are given respectively by
\begin{center}
\mbox{}
\!\!\!\!\!\!
\begin{tabular}{cc}
     \begin{tabular}{lll}
          $\Lzero  :=\textnormal{\tiny$\Lzeromatrix$}$, & $\Lone   :=\textnormal{\tiny$\Lonematrix$}$,    \\ \\
          $\Ltwo    :=\textnormal{\tiny$\Ltwomatrix$}$,  & $\Lthree :=\textnormal{\tiny$\Lthreematrix$}$,
     \end{tabular}
     $\Lfour:=\textnormal{\tiny$\Lfourmatrix$}$.
\end{tabular}
\end{center}



\subsection{The Vinberg diagrams}
\label{subsection:VinbergDiagrams}
\setcounter{theorem}{0}\setcounter{figure}{0}\setcounter{table}{0}

\begin{figure}
\begin{center}
\begin{tabular}{ccc}
$\Lzero$\;\;\begin{tabular}{c}\input{diagram-MAIN-VinbergLIV.tex} \end{tabular} &&
$\Lone$\;\;\begin{tabular}{c}\input{diagram-MAIN-VinbergLII.tex} \end{tabular}
\end{tabular}
\vskip 0.3cm
\begin{tabular}{ccc}
\begin{tabular}{l}
$\Lfour$\;\;\begin{tabular}{c}\input{diagram-MAIN-VinbergLVI.tex} \end{tabular} \\ \\ \\ \\
$\Ltwo$\;\;\begin{tabular}{c}\input{diagram-VinbergLI.tex}  \end{tabular}
\end{tabular}
&&
\mbox{}\quad\begin{tabular}{c}\input{diagram-VinbergLIII.tex}  \end{tabular}\;\;$\Lthree$
\end{tabular}
\end{center}
\caption{}
\label{diagram:VinbergDiagrams}
\vskip 0.3cm
\mbox{}
\end{figure}     
     
TheVinberg diagrams \cite{Vinberg73} of the reflection subgroups of the (integral) isometry groups $\P\Isom(\Lzero)$, $\ldots,$ $\P\Isom(\Lfour)$ are shown in Figure \ref{diagram:VinbergDiagrams}.
In these diagrams, the following convention is used: No bond between two nodes means the two corresponding hyperplanes meet orthogonally; a single bond means they meet with interior angle $\pi/3$; a double bond means the interior angle is $\pi/4$; a triple bond means the interior angle is $\pi/6$; a bond marked with $\infty$ means the two hyperplanes are parallel; a dotted bond means they are ultraparallel.  The number of subdivisions within each node is minus one-half of the squared norm of the corresponding root.  The labeling of the nodes of the diagrams for $\Ltwo$ and $\Lthree$ will be used in sections \ref{subsection:VinbergDiagramLII} and \ref{subsection:VinbergDiagramLIII}.  The common labeling of these two sets of nodes is for economy of notation; the two sets otherwise have no relation to each other.  Each of these five Vinberg diagrams has no symmetries, when norms of roots are taken into account.  This implies that each of $\P\Isom(\Lzero), \ldots, \P\Isom(\Lfour)$ is a discrete reflection subgroup of $\Isom(\RHfive)$.  Hence, Corollary \ref{NecessaryCondn1:ExistenceTypeII} applies to each of them.
              
Ignoring norms of roots, only the Vinberg diagrams of $\P\Isom(\Ltwo)$ and $\P\Isom(\Lthree)$ have a $(\Z/2\Z)$-symmetry, which implies (by Corollary \ref{NecessaryCondn1:ExistenceTypeII}) that
\begin{equation*}
       \Stab_{\P\Isom\,\Lambda}\left(\RHfive_{[\chi]}\right)
     =\P\Stab_{\Isom\,\Lambda}\left(\Fix\,\chi\right),
     \quad
     \textnormal{for} \; \chi = \chizero, \chione, \chifour.
\end{equation*}


\subsection{The invariants of the induced isometries on $V=\Lambda/(1+\i)\Lambda$}
\label{subsection:ChiZeroOneTwoThreeAreReal}
\setcounter{theorem}{0}\setcounter{figure}{0}\setcounter{table}{0}

Let $\phizero, \ldots, \phifour \in O(V,q)$ be the involutive isometries on $V = \Lambda/(1+\i)\Lambda$ induced by $\chizero,\ldots,\chifour$ respectively.  Then, straightforward computations show that the two invariants mentioned in Lemma \ref{LemmaTwoInvariants} of $\phizero$, $\ldots$, $\phifour$ are as tabulated in Table \ref{table:NormOneFixedVectors}.

\begin{table}
\begin{center}
\begin{tabular}{|c||c|c|}
\hline
                                  & $\dim_{\F_{2}}(\Fix(\,\cdot\,))$ & $\begin{array}{c}
                                                                                         \textnormal{\small number of norm-one} \\
                                                                                         \textnormal{\small vectors in $\Fix(\,\cdot\,)$}
                                                                                         \end{array}$ \\
                    \hline\hline
                    $\phizero$  & $6$ &  $\cdot$ \\ \hline
                    $\phione$   & $5$ &  $\cdot$ \\ \hline
                    $\phitwo$    & $4$ &  $8$      \\ \hline
                    $\phithree$ & $3$ &  $\cdot$ \\ \hline
                    $\phifour$   & $4$ &  $4$       \\ \hline
\end{tabular}
\end{center}
\caption{}
\label{table:NormOneFixedVectors}
\end{table}

\vskip 0.4cm
\begin{remark}\quad
Based on the computations presented so far, we may conclude that $\chizero, \ldots, \chithree$ correspond to real binary octics of types $0, \ldots, 3$ respectively.  It is also clear that $\chifour$ is induced by either real binary octics of type $4$ or antipodal binary octics. 
\end{remark}


\subsection{$\chifour$ is induced by real binary octics of type $4$}
\label{subsection:ChifourIsReal}
\setcounter{theorem}{0}\setcounter{figure}{0}\setcounter{table}{0}


We can determine that $\chifour$ is induced by real binary octics of type $4$ (rather than by antipodal binary octics) by the following observations:
\begin{itemize}
\item  Recall that the collection $\H$ of discriminant hyperplanes in $\CHfive$
          consists of orthogonal complements of vectors in $\Lambda$ of squared
          norm $-2$, and that the smooth points of $\H$ correspond to nodal
          binary octics, i.e. singular binary octics with one double point and no other
          singularities.
\item  One of the roots of $\Lfour$ is of the form $(1+\i)\,w$, where $w$ is a
          primitive vector in $\Lambda$ of squared norm $-2$.  The fundamental
          domain of $\P\Isom(\Lfour)$ therefore has one discriminant wall, and octics
          parametrized by $\RHfive_{[\chifour]}$ can deform to nodal ones.
\item  Antipodal octics can only deform to octics which are more singular than the
          nodal ones.  (See the proof of Lemma \ref{RealAntipIAAIClassesAreDisjoint}.)
\end{itemize}
It is now clear that $\chifour$ is induced by real binary octics of type $4$.


\begin{remark}\quad
As already mentioned in subsection \ref{subsection:VinbergDiagrams}, we have
\begin{equation*}
     \Stab_{\P\Isom\Lambda}(\RH_{[\chi]}) = \P\Stab_{\Isom\Lambda}(\Fix\,\chi),
     \quad\textnormal{for}\; \chi = \chizero, \chione, \chifour.
\end{equation*}
It remains to determine, for $\chi = \chitwo, \chithree$, whether
$\Stab_{\P\Isom\Lambda}(\RHfiveX)$ is equal to
$\P\Stab_{\Isom\Lambda}(\Fix\,\chi)$, or is isomorphic to
$\P\Stab_{\Isom\Lambda}(\Fix\,\chi) \rtimes (\Z/2\Z)$.
\end{remark}


\subsection{Comparing $\Stab_{\P\Isom\Lambda}(\RH_{[\chitwo]})$ and $\P\Stab_{\Isom\Lambda}(\Fix\,\chitwo)$}
\label{subsection:VinbergDiagramLII}
\setcounter{theorem}{0}\setcounter{figure}{0}\setcounter{table}{0}

Recall that the Vinberg diagram for $\P\Isom(\Fix\,\chitwo)$ is
\vskip 0.4cm
\begin{center} \input{diagram-VinbergLI.tex} \end{center}
\vskip 0.4cm
\noindent
where the roots $r_{1}, r_{2}, \ldots, r_{7}$ are labeled according to order of appearance in the Vinberg Algorithm.  The above diagram has only one symmetry (ignoring norms of roots): it is the $(\Z/2\Z)$-symmetry determined by exchanging the following $1$-dimensional subspaces:
\begin{equation*}
  \Re \cdot r_{1} \longleftrightarrow \Re \cdot r_{1}, \quad 
  \Re \cdot r_{2} \longleftrightarrow \Re \cdot r_{6}, \quad 
  \Re \cdot r_{3} \longleftrightarrow \Re \cdot r_{4}, \quad 
  \Re \cdot r_{5} \longleftrightarrow \Re \cdot r_{7}.
\end{equation*}
Recall also that the natural identification map (induced by projectivizing over $\C$) from $\Fix(\i\,\chitwo)\otimes_{\Z}\Re$ back to $\Fix(\chitwo)\otimes_{\Z}\Re$ is given by scalar multiplication by $1-\i$.  Taking all the above observations into account, we see that the group $\Stab_{\P\Isom\Lambda}\left(\RHfive_{[\chitwo]}\right)$ admits elements of type II if and only if the following conditions define an element $T \in \Isom(\Lambda)$ such that $[\,T\,]$ is a type II element of $\Stab_{\P\Isom\Lambda}\left(\RHfive_{[\chitwo]}\right)$:
\begin{equation*}
  (1-\i)\,T(r_{1}) = \pm\sqrt{2}\,r_{1},
\end{equation*} \vskip -0.4cm
\begin{equation}\label{Tconditions}
\begin{array}{lcl}
  (1-\i)\,T(r_{2}) = \pm\sqrt{2}\,r_{6},                   & \mbox{} &(1-\i)\,T(r_{6}) = \pm\sqrt{2}\,r_{2},      \\
  (1-\i)\,T(r_{3}) = \pm\dfrac{1}{\sqrt{2}}\,r_{4}, & \mbox{} &(1-\i)\,T(r_{4}) = \pm2\,\sqrt{2}\,r_{3},  \\
  (1-\i)\,T(r_{5}) = \pm\sqrt{2}\,r_{7},                   & \mbox{} &(1-\i)\,T(r_{7}) = \pm\sqrt{2}\,\,r_{5},
\end{array}
\end{equation}
where the signs of the right-hand-sides must be either all positive or all negative. Either case leads to a contradiction, which shows that $\Stab_{\P\Isom\Lambda}\left(\RHfive_{[\chitwo]}\right)$ has no type II elements.  We derive the contradiction for only the first case, the other case being completely analogous. We now make the following:
\begin{center}
{\bf CLAIM:}  There exists no such $T \in \Isom(\Lambda)$.
\end{center}
When expressed in the ``standard'' basis of $\Lambda$, the roots $r_{1}, \ldots, r_{7}$ are given, respectively from left to right, by the column vectors of the following matrix:
\begin{equation*}
     \textnormal{\scriptsize$\ZIrootsI$}.
\end{equation*}
Note that $r_{7}  \; = \;  r_{2} + 2\,r_{3} - r_{4} - r_{5} + r_{6}$.  Hence, the last condition in \eqref{Tconditions} implies
\begin{eqnarray*}
  \sqrt{2}\,r_{5} \;\, = \;\, (1-\i)\,T(r_{7})
  &=& (1-\i)\,T( r_{2} + 2\,r_{3} - r_{4} - r_{5} + r_{6}) \\
  &=& \sqrt{2}\,r_{6} + \sqrt{2}\,r_{4} - 2\,\sqrt{2}\,r_{3} - \sqrt{2}\,r_{7} + \sqrt{2}\,r_{2},
\end{eqnarray*}
which yields this alternative expression for $r_{7}$:\; $r_{7} = r_{2} - 2\,r_{3} + r_{4} - r_{5} + r_{6}$. Comparing with the original expression for $r_{7}$ in terms of $r_{2}, \ldots, r_{6}$, we see that
\begin{equation*}
  r_{2} + 2\,r_{3} - r_{4} - r_{5} + r_{6} \, = \, r_{7} \, = \, r_{2} - 2\,r_{3} + r_{4} - r_{5} + r_{6}
  \quad\Longrightarrow\quad
  2\,r_{3} = r_{4},  
\end{equation*}
which is a contradiction, since $r_{3}$ and $r_{4}$ are linearly independent over $\Zi$, in particular, over $\Z$.  The \claim is proved.  By Corollary \ref{NecessaryCondn1:ExistenceTypeII}, we may now conclude that $\Stab_{\P\Isom\Lambda}\left(\RHfive_{[\chitwo]}\right)$ has no elements of type II, and $\Stab_{\P\Isom\Lambda}\left(\RHfive_{[\chitwo]}\right) \, = \, \P\Stab_{\Isom\Lambda}\left(\Fix\,\chitwo\right)$.


\subsection{Comparing $\Stab_{\P\Isom\Lambda}(\RH_{[\chithree]})$ and $\P\Stab_{\Isom\Lambda}(\Fix\,\chithree)$}
\label{subsection:VinbergDiagramLIII}
\setcounter{theorem}{0}\setcounter{figure}{0}\setcounter{table}{0}

Recall the Vinberg diagram in this case from Figure \ref{diagram:VinbergDiagrams}.  Again, the roots $r_{1}, r_{2}, \ldots, r_{8}$ are labeled according to order of appearance in the Vinberg Algorithm.  In terms of the ``standard" basis for $\Lambda$, these roots are given, respectively from left to right, by the column vectors of the following matrix:
\begin{equation*}
      \textnormal{\scriptsize$\ZIrootsIII$}
\end{equation*}
The only symmetry (ignoring norms of roots) here is the $(\Z/2\Z)$-symmetry determined by exchanging the following $1$-dimensional subspaces:
\begin{equation*}
  \Re \cdot r_{1} \longleftrightarrow \Re \cdot r_{4}, \quad 
  \Re \cdot r_{2} \longleftrightarrow \Re \cdot r_{5}, \quad 
  \Re \cdot r_{3} \longleftrightarrow \Re \cdot r_{6}, \quad 
  \Re \cdot r_{7} \longleftrightarrow \Re \cdot r_{8}.
\end{equation*}
Therefore, $\Stab_{\P\Isom\Lambda}(\RH_{[\chithree]})$ has elements of type II if and only if the following conditions define an element $T \in \Isom(\Lambda)$ such that $[\,T\,]$ is an element of $\Stab_{\P\Isom\Lambda}(\RH_{[\chithree]})$ of type II:
{\small
\begin{equation*}
\begin{array}{llll}
(1-\i)\,T(r_{1}) = r_{4},     & (1-\i)\,T(r_{2}) = r_{5},      &   (1-\i)\,T(r_{3}) = r_{6},    &    (1-\i)\,T(r_{7}) = r_{8},  \\
(1-\i)\,T(r_{4}) = 2\,r_{1}, & (1-\i)\,T(r_{5}) = 2\,r_{2},  & (1-\i)\,T(r_{6}) = 2\,r_{3},  & (1-\i)\,T(r_{8}) = 2\,r_{7}.
\end{array}
\end{equation*}}

\noindent
Straightforward calculations now show that the above (overdetermined) set of conditions indeed defines such a $T \in \Isom(\Lambda)$ and we conclude that
\begin{equation*}
                   \Stab_{\P\Isom\Lambda}(\RH_{[\chithree]})
     \;=\;        \P\Stab_{\Isom\Lambda}(\Fix\,\chithree) \rtimes \langle [\,T\,] \rangle
     \;\cong\; \P\Stab_{\Isom\Lambda}(\Fix\,\chithree) \rtimes \left( \Z/2\Z \right).
\end{equation*}


%% file: diagram-MAIN-VinbergLIV.tex
\begin{picture}(130,40)

\put(  0,30){\circle{10}} \put(  0,25){\line(0,1){10}}
\put( 30,30){\circle{10}} \put( 30,25){\line(0,1){10}}
\put( 60,30){\circle{10}} \put( 60,25){\line(0,1){10}}
\put( 90,30){\circle{10}} \put( 90,25){\line(0,1){10}}
\put(120,30){\circle{10}} 

\put( 60, 0){\circle{10}} \put( 60,-5){\line(0,1){10}}


\put(  5,30){\line(1,0){20}}
\put( 35,30){\line(1,0){20}} 
\put( 65,30){\line(1,0){20}}
\put( 95,29){\line(1,0){20}}
\put( 95,31){\line(1,0){20}}

\put( 60,5){\line(0,1){20}}

\end{picture}

%% file: diagram-MAIN-VinbergLII.tex
\begin{picture}(160,40)

\put(  0,30){\circle{10}} \put(  0,25){\line(0,1){10}}
\put( 30,30){\circle{10}} \put( 30,25){\line(0,1){10}}
\put( 60,30){\circle{10}} \put( 60,25){\line(0,1){10}}
\put( 90,30){\circle{10}} \put( 90,25){\line(0,1){10}}
\put(120,30){\circle{10}} \put(120,25){\line(0,1){10}} 
\put(150,30){\circle{10}} 

\put(120, 0){\circle{10}} \put(120,-5){\line(0,1){10}}



\put(  5,30){\line(1,0){20}} \put(10,32){$\infty$}
\put( 35,30){\line(1,0){20}} 
\put( 65,30){\line(1,0){20}}
\put( 95,30){\line(1,0){20}}

\put(125,29){\line(1,0){20}}
\put(125,31){\line(1,0){20}}

\put(120,5){\line(0,1){20}}

\end{picture}

%% file: diagram-MAIN-VinbergLVI.tex
\begin{picture}(160,10)

\put(  0,0){\circle{10}} \put(  0,-5){\line(0,1){10}}

\thinlines
\put( 30,0){\circle{10}} \put( 30,-5){\line(0,1){10}}
\put( 60,0){\circle{10}} \put( 60,-5){\line(0,1){10}}
\put( 90,0){\circle{10}} \put( 90,-5){\line(0,1){10}}
\put(120,0){\circle{10}} \put(120,-5){\line(0,1){10}} \put(115,0){\line(1,0){10}}
\put(150,0){\circle{10}} \put(150,-5){\line(0,1){10}} \put(145,0){\line(1,0){10}}

\put(  5, 0){\line(1,0){20}}
\put( 35, 0){\line(1,0){20}}
\put( 65, 0){\line(1,0){20}}
\put( 95,-1){\line(1,0){20}}
\put( 95, 1){\line(1,0){20}}
\put(125, 0){\line(1,0){20}}

\end{picture}

%% file: diagram-VinbergLI.tex
\begin{picture}(190,10)

\put(  0,0){\circle{10}} \put(  0,-5){\line(0,1){10}}
\put( 30,0){\circle{10}} 
\put( 60,0){\circle{10}} \put( 60,-5){\line(0,1){10}}
\put( 90,0){\circle{10}} \put( 90,-5){\line(0,1){10}}
\put(120,0){\circle{10}} \put(120,-5){\line(0,1){10}} 
\put(150,0){\circle{10}} \put(150,-5){\line(0,1){10}} \put(145,0){\line(1,0){10}}
\put(180,0){\circle{10}} \put(180,-5){\line(0,1){10}} 


\put( 87,10){\scriptsize$r_{1}$}
\put( 57,10){\scriptsize$r_{2}$}
\put( 27,10){\scriptsize$r_{3}$}
\put(117,10){\scriptsize$r_{6}$}
\put(147,10){\scriptsize$r_{4}$}
\put(177,10){\scriptsize$r_{5}$}
\put( -3,10){\scriptsize$r_{7}$}


\put(  5,-1){\line(1,0){20}}
\put(  5, 1){\line(1,0){20}}

\put( 35,-1){\line(1,0){20}}
\put( 35, 1){\line(1,0){20}}

\put( 65, 0){\line(1,0){20}}
\put( 95, 0){\line(1,0){20}}

\put(125,-1){\line(1,0){20}}
\put(125, 1){\line(1,0){20}}

\put(155,-1){\line(1,0){20}}
\put(155, 1){\line(1,0){20}}

\end{picture}

%% file: diagram-VinbergLIII.tex
\begin{picture}(100,100)

\put(  0,30){\circle{10}} 
\put( 30, 0){\circle{10}} 
\put( 60, 0){\circle{10}} \put( 60,-5){\line(0,1){10}}
\put( 90,30){\circle{10}} \put( 90,25){\line(0,1){10}}

\put(  0,60){\circle{10}} \put(  0,55){\line(0,1){10}}
\put( 30,90){\circle{10}} \put( 30,85){\line(0,1){10}}
\put( 60,90){\circle{10}} \put( 60,85){\line(0,1){10}} \put(55,90){\line(1,0){10}}
\put( 90,60){\circle{10}} \put( 90,55){\line(0,1){10}} \put(85,60){\line(1,0){10}}

\put(-14, 28){\textnormal{\scriptsize$r_{2}$}} 
\put(-14, 58){\textnormal{\scriptsize$r_{1}$}}
\put( 27, 98){\textnormal{\scriptsize$r_{6}$}}
\put( 57, 98){\textnormal{\scriptsize$r_{3}$}}
\put( 97, 58){\textnormal{\scriptsize$r_{4}$}}
\put( 97, 28){\textnormal{\scriptsize$r_{5}$}}
\put( 27,-12){\textnormal{\scriptsize$r_{7}$}}
\put( 57,-12){\textnormal{\scriptsize$r_{8}$}}


\put( 35,89){\line(1,0){20}}
\put( 35,91){\line(1,0){20}}

\put(-1,35){\line(0,1){20}}
\put( 1,35){\line(0,1){20}}

\put(89,35){\line(0,1){20}}
\put(91,35){\line(0,1){20}}

\put(3, 64){\line( 1,1){23}}
\put(63, 4){\line( 1,1){23}} \put(75,5){$\infty$}

\put(87,64){\line(-1,1){23}} 
\put(27, 4){\line(-1,1){23}} \put(5,5){$\infty$}

\thicklines\dottedline{3}(36,0)(54,0)\thinlines

\end{picture}

%% file: NonOrbifoldPoints-nocurves.tex
\section{$\MsR$ Is Not a Real Hyperbolic Orbifold}
\label{section:nonhyperbolicity}
\setcounter{theorem}{0}\setcounter{figure}{0}\setcounter{table}{0}

A (singular) complex binary octic is said to be \emph{cuspidal} if it has exactly one triple point and no other singularities.  Note that the antipodal map on $\CPone$ cannot preserve any cuspidal octic, and smooth antipodal octics cannot deform to a cuspidal octic.

In this section, we show that $\MsR$ is not a hyperbolic orbifold by proving that its metric-space structure around moduli points in the stratum $\Delta^{0,1}_{\Re}$ cannot be described as the quotient of (an open subset of) a Riemannian manifold by the action of a finite group of isometries, where $\Delta^{0,1}_{\Re}$ denotes the stratum of moduli points that corresponds to real cuspidal octics.

\subsection{The vanishing $(\sigma^{2}=-1)$-cohomology of $X_{p}$ for a cuspidal octic $p$}
\setcounter{theorem}{0}\setcounter{figure}{0}\setcounter{table}{0}

\begin{lemma}\label{CuspidalSummand}\quad
The vanishing $(\sigma^{2}=-1)$-cohomology $\Lz(p)$ corresponding to a real cuspidal octic $p$ is an orthogonal summand of $\Lambda(X_{p_{0}}) \cong \Lambda$, where $p_{0}$ is any reference smooth complex binary octic, and $\Lz(p)$ is isometric to the following $\Zi$-lattice of $\Zi$-rank two:
\begin{equation*}
  \Lz := \left(\; \Zi^{2} \; , \; \textnormal{\scriptsize$\cuspsummand$} \;\right).
\end{equation*}
\end{lemma}
\proofoutline
We can describe locally the deformation from a smooth real binary octic to a singular real cuspidal octic by examining:
\begin{equation*}
   p_{a_{0},a_{1}}(x) \,=\, (x^{3} + a_{1}x + a_{0})\cdot r(x),
\end{equation*}
as $a_{0}, a_{1} \rightarrow 0$, where $a_{0}, a_{1}\in\Re$, and $r(x)$ is a polynomial in $x$ of degree five with real coefficients and has no common roots with $x^{3} + a_{1}x + a_{0}$.

Recall that our Hodge-theory set-up arises from the $(\sigma^{2}=-1)$-eigenspace of the cyclic action on the
cohomology of the four-fold cyclic cover $X_{p}$ of $\CPone$ branched over the roots of any smooth binary octic $p$.  The intersection form of the vanishing $(\sigma^{2}=-1)$-homology of $X_{p}$ corresponding to the singularity of one triple point and no other singularities can be described locally by the corresponding  vanishing $(\sigma^{2}=-1)$-homology of the singularity of $p_{0,0}(x)$ above.  A simple pictorial argument shows that this $(\sigma^{2}=-1)$-homology indeed has the $\Zi$-lattice structure of $\Lz$. \qed

\begin{remark}\quad
A combinatorial argument shows that $\Isom(\Lz)$ has $96$ elements; see \cite{Chu:PhDthesis}.  In fact, $\Isom(\Lz)$ is isomorphic to the group $\Bthreefour$, the group obtained from the braid group $B_{3}$ on three strands by imposing an order-four condition on the ``standard" generators.  This fact is intuitively clear since $\Lz$ is the intersection form of the vanishing (co)homology (strictly speaking, a certain eigensubspace of it) of an order-four branched cyclic cover of $\CPone$ corresponding to a singularity of the coalescence of three branch points.
\end{remark}


\subsection{The local quotient structure near the period of a generic cuspidal octic}
\label{subsection:CuspidalLocalStructure}
\setcounter{theorem}{0}\setcounter{figure}{0}\setcounter{table}{0}

Let $x \in \CHfive$ be the period of a generic cuspidal real binary octic $p$.  Since a triple point can be thought of as the ``limit'' of two nodes, it is intuitively clear that the vanishing $(\sigma^{2}=-1)$-cohomology $\Lz(p)$ corresponding to the singularity of $p$ has $\Zi$-rank two, and that the corresponding local monodromy group should have a natural representation on $\Zi^{2}$.  Proposition \ref{CuspidalSummand} asserts that $\Lz(p)$ is indeed isometric to the abstract $\Zi$-lattice $\Lz$ and the corresponding local monodromy group is isomorphic to $\Isom\left(\Lz\right)$.

Since $p$ is real and singular, its period $x \in \CHfive$ lies on a collection of (more than one) integral copies of $\RHfive$.  The common intersection of this collection of integral copies of $\RHfive$ is a totally real copy of $\RH^{3}$ in $\CHfive$.  
Since cuspidal octics are stable, we know from Geometric Invariant Theory that $\Stab_{\P\Isom\,\Lambda}(x)$ is a finite group.  Being isometries of $\Lambda$, elements of $\Stab_{\P\Isom(\Lambda)}(x)$ preserve both $\Lz(p)$ and $\Lz(p)^{\perp}$ individually.  Similarly, involutive anti-isometries of $\Lambda$ that preserve $\Lz(p)$ must also preserve $\Lz(p)^{\perp}$.  Each such anti-isometry of $\Lambda$ of course restricts to an involutive anti-isometry on $\Lz(p)$ and $\Lz(p)^{\perp}$ individually.  Furthermore, genericity of $p$ implies that $\Stab_{\P\Isom\,\Lambda}(x)$ acts trivially on $\Lz(p)^{\perp}$.  

All of the above implies that the real period $x$ has a neighborhood $\mathcal{U}$ in $ \Ks$ which is homeomorphic to an open neighborhood of $(\textnormal{origin},\textnormal{generic point})$ in
\begin{equation}\label{eqn:LocalBeforeQuotientStructure}
\left\{\begin{array}{c}
\textnormal{\small the union of the fix-point-sets of}\\
\textnormal{\small all involutive anti-isometries of $\Lz$}
\end{array}\right\}
\;\textnormal{\huge$\times$}\;
\RH^{3}.
\end{equation}
We emphasize however that the natural metric-space structure on $\mathcal{U}$ (i.e. the one inherited from that of $\Ks$) is not isometric to the product metric space of the two factors in \eqref{eqn:LocalBeforeQuotientStructure}.  The action of $\Stab_{\P\Isom\,\Lambda}(x)$ on $\mathcal{U}$ gives rise to a local quotient homeomorphic to an open neighborhood of the image of $(\textnormal{origin},\textnormal{generic point})$ in
\begin{equation}\label{eqn:LocalQuotientStructure}
\left(
\Isom(\Lz)\left\backslash
\left\{\!\begin{array}{c}
\textnormal{\small the union of the fix-point-sets of}\\
\textnormal{\small all involutive anti-isometries of $\Lz$}
\end{array}\!\right\}
\right.
\right)
\textnormal{\huge$\times$}\;
\RH^{3}.
\end{equation}
We may choose coordinates on $\CPone$ so that the triple point of the real cuspidal octic $p(z)$ occurs at $0 \in \CPone = \C\cup\{\infty\}$.  In non-homogeneous coordinates, we may write $p(z) = z^{3} \cdot r(z)$, where $r(z)$ is a polynomial in $z$ of degree $5$ with distinct roots, each distinct from $0 \in \CPone$.  Then the first factor in \eqref{eqn:LocalQuotientStructure} describes the desingularization of the triple point into smooth real 3-point configurations, whereas the second factor describes the deformation of the roots of $r(z)$.

The first factor in \eqref{eqn:LocalQuotientStructure} turns out to be a flat two-real-dimensional cone, obtained by gluing together two flat two-dimensional wedges. Similarly, the subspace of $\mathcal{U}$ that topologically corresponds to this first factor can be abstractly described by gluing together two real-hyperbolic two-dimensional wedges.  We may therefore consider the cone angles at the vertices of these two cones (the former is flat; the latter is not).   Noting that $T_{x}\CHfive$ is isometric to the orthogonal complement of $x$ in $\C^{1,5}$, we see that these two cone angles are equal.  We will show in the following subsections that the common cone angle value is $3\pi/4$.  

The $3\pi/4$ cone angle shows that $x$ cannot be a real-Riemannian orbifold point of $\MsR$, since otherwise the cone angle would have to be an integral submultiple of $\pi$.  This in particular implies that the Allcock-Carlson-Toledo construction $\MsR$ of the moduli space of stable real binary octics cannot be a real-hyperbolic orbifold, in contrast to the cases of real cubic surfaces \cite{ACT:realcubics} and real sextics \cite{ACT:realsextics}.



\subsection{The two isometry classes of involutive anti-isometries of $\Lz$}
\setcounter{theorem}{0}\setcounter{figure}{0}\setcounter{table}{0}

We now begin the computation of the cone angle of the first factor in \eqref{eqn:LocalQuotientStructure}.

\begin{proposition}\label{IAAI:Lz}\quad
$\Lz$ admits exactly two $\Isom(\Lz)$-conjugacy classes of involutive anti-isometries, represented by:
\begin{equation*}
     \kone  \left(\xvector\right) = \konematrix,
     \quad\textnormal{and}\quad 
     \kthree\left(\xvector\right) = \kthreematrix.
\end{equation*}
\end{proposition}
\proofoutline
\begin{enumerate}

\item  \textit{The maps $\kone : \Lz\longrightarrow\Lz$ and $\kthree : \Lz\longrightarrow\Lz$ are involutive anti-isometries of $\Lz$.}  This can be verified with straightforward calculations.

\item  \textit{There are exactly 36 involutive anti-isometries of $\Lz$.}  First note that every involutive anti-isometry $\kappa$ of $\Lz$ is of the form $A \circ \kone$, for some isometry $A \in \Isom\,\Lz$.  This is because $\kappa\circ\kone \in \Isom\,\Lz$, for every involutive anti-isometry $\kappa$ of $\Lz$.  On the other hand, note that involutiveness of $\kappa$ implies $\textnormal{id}_{\Lz} = \kappa^{2} = \left( A \circ \kone \right) \circ \left( A \circ \kone \right)$; equivalently, $\kone \; = \; A \circ \kone \circ A$.  We therefore have
\begin{equation*}
  \IAAI(\Lz) \; := \;
  \left\{\begin{array}{c}
    \textnormal{\small involutive} \\ \textnormal{\small anti-isometries} \\ \textnormal{\small of $\Lz$}
  \end{array}\right\}
  \; = \;
  \left\{\;
     A \circ \kone \;
     \left\vert\;
     \begin{array}{c}A \in \Isom\,\Lz, \\ \kone \; = \; A \circ \kone \circ A\end{array}
     \;\right.
  \;\right\}.
\end{equation*}
Hence, $|\,\IAAI(\Lz)\,| \leq |\,\Isom\,\Lz\,| = 96$.  In fact, calculations show that there exactly $36$ isometries $A \in \Isom\,\Lz$ which satisfy $\kone = A \circ \kone \circ A$.  Hence, $|\,\IAAI(\Lz)\,| = 36$.

\item  \textit{There are exactly two $(\Isom\,\Lz)$-classes of involutive anti-isometries of $\Lz$, represented by $\kone$ and $\kthree$ respectively.}  By Lemmas \ref{Fixkone} and \ref{Fixkthree}, $\Fix\,\kone$ and $\Fix\,\kthree$ are not isometric, we know there are at least two $(\Isom\,\Lz)$-classes of involutive anti-isometries of $\Lz$.  On the other hand, calculations show that each of the 36 involutive anti-isometries of $\Lz$ is $(\Isom\,\Lz)$-conjugate to either $\kone$ or $\kthree$.  Hence, we conclude that there are exactly two $(\Isom\,\Lz)$-classes.  \qed

\end{enumerate}

\begin{lemma}\label{Fixkone}\quad
   $\Fix(\kone)$ is isometric to $\scriptsize\left(\; \Z\oplus\Z \,,\, \fixkonemetric \;\right)$
   and the action
   \begin{equation*} \Stab_{\Isom\,\Lz}(\Fix\,\kone) \acts \Fix\,\kone \end{equation*}
   is equivalent to the action $(\Z/2\times\Z/2)\,\acts\,\Z^{2}$.  Consequently,
   {\large$_{\Stab_{\Isom\,\Lz}(\Fix\,\kone)} \backslash \left(\Fix(\kone)\otimes_{\Z}\Re\right)$}
   is isometric to the ``\,$90^{\circ}$-wedge''
   \begin{equation*}
      \left\{\; (x,y)\in\Re^{2}  \;\left|\; x, y \geq 0 \right.\;\right\},
   \end{equation*}
   where the latter has the usual Euclidean metric.
\end{lemma}
\proofoutline It can be shown that $\uone$ $=$ {\scriptsize$\left(\begin{array}{r} -1-\i \\ 1-\i \end{array}\right)$}, and $\utwo$ $=$ {\scriptsize$\left(\begin{array}{r} \i \\ \i \end{array}\right)$} form a $\Z$-basis for $\Fix(\kone)$.  Straightforward calculations show that their inner product matrix is $\scriptsize\fixkonemetric$. Consequently,
\begin{equation*}
   \Fix(\kone) \;\; = \;\;
   \left(\;
     \Z\cdot\uone\,\oplus\,\Z\cdot\utwo \;,\;
     \textnormal{\scriptsize$\fixkonemetric$}
   \;\right)
   \;\;\cong\;\;
   \left(\;
     \Z\oplus\Z \;,\;
     \textnormal{\scriptsize$\fixkonemetric$}
   \;\right).
\end{equation*}
The rest of the Lemma follows easily from the fact that the intersection form on $\Fix(\kone)$ is $\diag(-4,-2)$. \qed

\begin{lemma}\label{Fixkthree}\quad
$\Fix(\kthree)$ is isometric to $\scriptsize\left(\; \Z\oplus\Z \,,\, \fixkthreemetric \;\right)$ and the action
\begin{equation*} \Stab_{\Isom\,\Lz}(\Fix\,\kthree) \acts \Fix\,\kthree \end{equation*}
is equivalent to the action $D_{4}\,\acts\,\Z^{2}$, where $D_{4}$ is the dihedral group of eight elements.  Consequently, {\large$_{\Stab_{\Isom\,\Lz}(\Fix\,\kthree)} \backslash \left(\Fix(\kthree)\otimes_{\Z}\Re\right)$}
is isometric to the ``$45^{\circ}$-wedge''
\begin{equation*}
     \left\{\; (x,y)\in\Re^{2}  \;\left|\; 0 \leq y \leq x \right.\;\right\},
\end{equation*}
where the latter has the usual Euclidean metric.
\end{lemma}
\proofoutline It can shown that $\vone$ $=$ {\scriptsize$\left(\begin{array}{c} 1+\i \\ \i \end{array}\right)$}, and $\vtwo$ $=$ {\scriptsize$\left(\begin{array}{c} 0 \\ -\,\i \end{array}\right)$} form a $\Z$-basis for $\Fix(\kthree)$ and
that their inner product matrix is $\scriptsize\fixkthreemetric$.  This immediately shows that
\begin{equation*}
   \Fix(\kthree) \;\; = \;\;
   \left(\;
     \Z\cdot\vone\,\oplus\,\Z\cdot\vtwo \;,\;
     \textnormal{\scriptsize$\fixkthreemetric$}
   \;\right)
   \;\;\cong\;\;
   \left(\;
     \Z\oplus\Z \;,\;
     \textnormal{\scriptsize$\fixkthreemetric$}
   \;\right).
\end{equation*}
The rest of the Lemma is clear since the quadratic form on $\Fix(\kthree)$ is $\diag(-2,-2)$.  \qed


\subsection{The gluing of the fixed-point-sets of the involutive anti-isometries of $\Lz$ induced by the action of $\Isom(\Lz)$}
\setcounter{theorem}{0}\setcounter{figure}{0}\setcounter{table}{0}

Recall that $\Fix(\kone) = \Z\cdot\uone\,\oplus\,\Z\cdot\utwo$ and $\Fix(\kthree) = \Z\cdot\vone\,\oplus\,\Z\cdot\vtwo$,  where $\uone$ $=$ {\scriptsize$\left(\begin{array}{r} -1-\i \\ 1-\i \end{array}\right)$}, $\utwo$ $=$ {\scriptsize$\left(\begin{array}{r} \i \\ \i \end{array}\right)$}, $\vone$ $=$ {\scriptsize$\left(\begin{array}{c} 1+\i \\ \i \end{array}\right)$}, and $\vtwo$ $=$ {\scriptsize$\left(\begin{array}{c} 0 \\ -\,\i \end{array}\right)$}.  Recall also that\, {\large$_{\Stab_{\Isom\,\Lz}(\Fix\,\kone)}\,\backslash\,(\Fix\,\kone\otimes\Re)$}\, is a $90^{\circ}$-wedge, whereas\, {\large$_{\Stab_{\Isom\,\Lz}(\Fix\,\kthree)}\,\backslash\,(\Fix\,\kthree\otimes\Re)$}\, is a $45^{\circ}$-wedge.  Now, define
\begin{equation*}
   \vthree \; := \; \vone + \vtwo \; = \; \left(\begin{array}{c} 1+\i \\ 0 \end{array}\right)
   \; \in \; \Fix(\kthree).
\end{equation*}

\begin{lemma}\label{IsomLz:gluing}\mbox{}
\begin{enumerate}
   \item  $\uone$ and $\vthree$ belong to the same $(\Isom\,\Lz)$-orbit, and
   \item  $\utwo$ and $\vtwo$   belong to the same $(\Isom\,\Lz)$-orbit.
\end{enumerate}
\end{lemma}
\proof  Let $A_{1}$ $:=$ {\scriptsize$\left[\begin{array}{rrr} -1 &&  0 \\ -\i && \i \end{array}\right]$} and $A_{2}$ $:=$ {\scriptsize$\left[\begin{array}{rrr}  0 && -1 \\   1 && -1 \end{array}\right]$}.  Straightforward calculations show that $A_{1}\cdot\vthree = \uone$, $A_{2}\cdot\vtwo = \utwo$, and $A_{1}, A_{2} \in \Isom\,\Lz$.  This completes the proof.  \qed

\vskip 0.5cm
We are now ready to give an explicit description of the local quotient mentioned in section \ref{subsection:CuspidalLocalStructure}:
\begin{proposition}\label{Propn:LocalCuspidalStructure}\quad
A generic point on the stratum $\Delta^{0,1}_{\Re}$ of the Allcock-Carlson-Toledo construction $\MsR$ of the moduli space of stable real binary octics has a neighborhood which is homeomorphic to an open neighborhood of $(\textnormal{vertex},\textnormal{generic point})$ in $\mathcal{C}\times\RH^{3}$, where $\mathcal{C}$ is a real-two-dimensional cone.  Furthermore, the induced metric-space structure on $\mathcal{C}\times\{\textnormal{generic point}\}$ is such that the cone angle at the vertex is $3\pi/4$.
\end{proposition}

\proofoutline  The local topological product structure is given by \eqref{eqn:LocalQuotientStructure}.  Lemmas \ref{Fixkone}, \ref{Fixkthree}, and \ref{IsomLz:gluing} together imply that the first factor in \eqref{eqn:LocalQuotientStructure}
is isometric to an open neighborhood of the vertex in the real-two-dimensional cone obtained from gluing a Euclidean $45^{\circ}$-wedge with a Euclidean $90^{\circ}$-wedge along the edges as shown in Figure \ref{diagram:TwoWedges}.
\begin{figure}
\begin{center}
\input{diagram-wedges-nocurves.tex} \vskip -1.5cm
\caption{\small\it The gluing of the two wedges is given by identifying $\utwo$ with $\vtwo$, as well as $\uone$ with $\vthree$.  This gives rise to a $135^{\circ}$-wedge.}
\label{diagram:TwoWedges}
\end{center}
\end{figure}
We remark that the only non-manifold point in the above construction is the ``vertex"; in particular, the points along the edges spanned by $\uone$ (or $\vthree$) and $\utwo$ (or $\vtwo$) are manifold points (except the vertex itself).

We now see the cone angle at the vertex of the non-flat cone $\mathcal{C}\times\{\textnormal{generic point}\}$ is $3\pi/4$ from the fact that the flat cone constructed in the preceding paragraph is the infinitesimal representation of $\mathcal{C}\times\{\textnormal{generic point}\}$ at the vertex.  \qed

\begin{remark}\quad
The local metric-space structure of $\left.\Stab_{\P\Isom\,\Lambda}(x)\,\right\backslash\,\mathcal{U}$ mentioned in subsection \ref{subsection:CuspidalLocalStructure} is in fact a fibration over an open set of $\RH^{3}$ with fiber being a two-dimensional real-hyperbolic cone with cone angle $3\pi/4$.  Figure \ref{diagram:TwoWedges} can also be used to help illustrate this as follows:  The $45^{\circ}$-wedge on the left depicts a wedge in $\RHfive$ bounded by two copies of $\RH^{4}$ which intersect with interior angle $\pi/4$ along a common copy of $\RH^{3}$, which is represented by the vertex.  Each of the two boundary copies of $\RH^{4}$ is represented by the vertex together with either the vector $\vtwo$ or $\vthree$.  Similarly, the $90^{\circ}$-wedge on the right also depicts a wedge in $\RHfive$, except that the interior angle there is $\pi/2$.  These two wedges of $\RHfive$ are glued together by identifying the copies of $\RH^{4}$ represented by $\utwo$ and $\vtwo$, as well as by identifying those represented by $\uone$ and $\vthree$.  The two copies of $\RH^{3}$ represented by the two vertices actually coincide in $\Ks$.  This copy of $\RH^{3}$ (the base of the fibration) parametrizes real configurations of 1 triple point and 5 single points.  The interior of the wedge on the left parametrizes the resolutions of the real triple point into three real points (with the remaining generic real $5$-point configuration held fixed), whereas that on the right parametrizes the resolutions of the real triple point into configurations of one real point and one complex conjugate pair. 
\end{remark}

\begin{corollary}\label{Corollary:MsRIsNOTHyperbolic}\quad
$\MsR$ is not a Riemannian orbifold; in particular, it cannot be a real hyperbolic orbifold.
\end{corollary}
\proof  Simply note that the local angle of any two-dimensional orbifold must be an integral submultiple of $\pi$, but $3\pi/4$ is not an integral submultiple of $\pi$.  \qed

%% file: diagram-wedges-nocurves.tex

\begin{picture}(400,200)

\put(235,165){$_{\Stab_{\Isom\,\Lz}(\Fix\,\kone)}\backslash\left(\Fix\,\kone\,\otimes\Re\right)$}

{\thicklines\put(250,45){\vector(1,0){71}}}
\put(315,35){$\uone$}

{\thicklines\put(250,45){\vector(0,1){50}}}
\put(235,90){$\utwo$}


\put(33,165){$_{\Stab_{\Isom\,\Lz}(\Fix\,\kthree)}\backslash\left(\Fix\,\kthree\,\otimes\Re\right)$}

{\thicklines\put(150,45){\vector(-1,0){50}}}
\put(95,35){$\vone$}

{\thicklines\put(150,45){\vector(0,1){50}}}
\put(155,90){$\vtwo$}

{\thicklines\put(150,45){\vector(-1,1){48}}}
\put(40,85){$\vthree = \vone + \vtwo$}


\qbezier(170,95)(200,125)(230,95)
\put(175,100){\vector(-1,-1){5}}
\put(225,100){\vector(1,-1){5}}

\qbezier(95,80)(-10,-20)(305,35)
\put(90,75){\vector(1,1){5}}
\put(302,34){\vector(3,1){5}}








\dottedline{2}(250, 45)(350,145)

\dottedline{2}(250, 47)(348,145)
\dottedline{2}(250, 49)(346,145)
\dottedline{2}(250, 51)(344,145)
\dottedline{2}(250, 53)(342,145)
\dottedline{2}(250, 55)(340,145)

\dottedline{2}(250, 57)(338,145)
\dottedline{2}(250, 59)(336,145)
\dottedline{2}(250, 61)(334,145)
\dottedline{2}(250, 63)(332,145)
\dottedline{2}(250, 65)(330,145)

\dottedline{2}(250, 67)(328,145)
\dottedline{2}(250, 69)(326,145)
\dottedline{2}(250, 71)(324,145)
\dottedline{2}(250, 73)(322,145)
\dottedline{2}(250, 75)(320,145)

\dottedline{2}(250, 77)(318,145)
\dottedline{2}(250, 79)(316,145)
\dottedline{2}(250, 81)(314,145)
\dottedline{2}(250, 83)(312,145)
\dottedline{2}(250, 85)(310,145)

\dottedline{2}(250, 87)(308,145)
\dottedline{2}(250, 89)(306,145)
\dottedline{2}(250, 91)(304,145)
\dottedline{2}(250, 93)(302,145)
\dottedline{2}(250, 95)(300,145)

\dottedline{2}(250, 97)(298,145)
\dottedline{2}(250, 99)(296,145)
\dottedline{2}(250,101)(294,145)
\dottedline{2}(250,103)(292,145)
\dottedline{2}(250,105)(290,145)

\dottedline{2}(250,107)(288,145)
\dottedline{2}(250,109)(286,145)
\dottedline{2}(250,111)(284,145)
\dottedline{2}(250,113)(282,145)
\dottedline{2}(250,115)(280,145)

\dottedline{2}(250,117)(278,145)
\dottedline{2}(250,119)(276,145)
\dottedline{2}(250,121)(274,145)
\dottedline{2}(250,123)(272,145)
\dottedline{2}(250,125)(270,145)

\dottedline{2}(250,127)(268,145)
\dottedline{2}(250,129)(266,145)
\dottedline{2}(250,131)(264,145)
\dottedline{2}(250,133)(262,145)
\dottedline{2}(250,135)(260,145)

\dottedline{2}(250,137)(258,145)
\dottedline{2}(250,139)(256,145)
\dottedline{2}(250,141)(254,145)
\dottedline{2}(250,143)(252,145)
\dottedline{2}(250,145)(250,145)


\dottedline{2}(250,45)(350,145)
\dottedline{2}(252,45)(350,143)
\dottedline{2}(254,45)(350,141)
\dottedline{2}(256,45)(350,139)
\dottedline{2}(258,45)(350,137)
\dottedline{2}(260,45)(350,135)

\dottedline{2}(262,45)(350,133)
\dottedline{2}(264,45)(350,131)
\dottedline{2}(266,45)(350,129)
\dottedline{2}(268,45)(350,127)
\dottedline{2}(270,45)(350,125)

\dottedline{2}(272,45)(350,123)
\dottedline{2}(274,45)(350,121)
\dottedline{2}(276,45)(350,119)
\dottedline{2}(278,45)(350,117)
\dottedline{2}(280,45)(350,115)

\dottedline{2}(282,45)(350,113)
\dottedline{2}(284,45)(350,111)
\dottedline{2}(286,45)(350,109)
\dottedline{2}(288,45)(350,107)
\dottedline{2}(290,45)(350,105)

\dottedline{2}(292,45)(350,103)
\dottedline{2}(294,45)(350,101)
\dottedline{2}(296,45)(350, 99)
\dottedline{2}(298,45)(350, 97)
\dottedline{2}(300,45)(350, 95)

\dottedline{2}(302,45)(350, 93)
\dottedline{2}(304,45)(350, 91)
\dottedline{2}(306,45)(350, 89)
\dottedline{2}(308,45)(350, 87)
\dottedline{2}(310,45)(350, 85)

\dottedline{2}(312,45)(350, 83)
\dottedline{2}(314,45)(350, 81)
\dottedline{2}(316,45)(350, 79)
\dottedline{2}(318,45)(350, 77)
\dottedline{2}(320,45)(350, 75)

\dottedline{2}(322,45)(350, 73)
\dottedline{2}(324,45)(350, 71)
\dottedline{2}(326,45)(350, 69)
\dottedline{2}(328,45)(350, 67)
\dottedline{2}(330,45)(350, 65)

\dottedline{2}(332,45)(350, 63)
\dottedline{2}(334,45)(350, 61)
\dottedline{2}(336,45)(350, 59)
\dottedline{2}(338,45)(350, 57)
\dottedline{2}(340,45)(350, 55)

\dottedline{2}(342,45)(350, 53)
\dottedline{2}(344,45)(350, 51)
\dottedline{2}(346,45)(350, 49)
\dottedline{2}(348,45)(350, 47)
\dottedline{2}(350,45)(350, 45)



\dottedline{4}(150,45)(50,145)
\dottedline{4}(150,47)(52,145)
\dottedline{4}(150,49)(54,145)
\dottedline{4}(150,51)(56,145)
\dottedline{4}(150,53)(58,145)

\dottedline{4}(150,55)(60,145)
\dottedline{4}(150,57)(62,145)
\dottedline{4}(150,59)(64,145)
\dottedline{4}(150,61)(66,145)
\dottedline{4}(150,63)(68,145)

\dottedline{4}(150,65)(70,145)
\dottedline{4}(150,67)(72,145)
\dottedline{4}(150,69)(74,145)
\dottedline{4}(150,71)(76,145)
\dottedline{4}(150,73)(78,145)

\dottedline{4}(150,75)(80,145)
\dottedline{4}(150,77)(82,145)
\dottedline{4}(150,79)(84,145)
\dottedline{4}(150,81)(86,145)
\dottedline{4}(150,83)(88,145)

\dottedline{4}(150,85)(90,145)
\dottedline{4}(150,87)(92,145)
\dottedline{4}(150,89)(94,145)
\dottedline{4}(150,91)(96,145)
\dottedline{4}(150,93)(98,145)

\dottedline{4}(150, 95)(100,145)
\dottedline{4}(150, 97)(102,145)
\dottedline{4}(150, 99)(104,145)
\dottedline{4}(150,101)(106,145)
\dottedline{4}(150,103)(108,145)

\dottedline{4}(150,105)(110,145)
\dottedline{4}(150,107)(112,145)
\dottedline{4}(150,109)(114,145)
\dottedline{4}(150,111)(116,145)
\dottedline{4}(150,113)(118,145)

\dottedline{4}(150,115)(120,145)
\dottedline{4}(150,117)(122,145)
\dottedline{4}(150,119)(124,145)
\dottedline{4}(150,121)(126,145)
\dottedline{4}(150,123)(128,145)

\dottedline{4}(150,125)(130,145)
\dottedline{4}(150,127)(132,145)
\dottedline{4}(150,129)(134,145)
\dottedline{4}(150,131)(136,145)
\dottedline{4}(150,133)(138,145)

\dottedline{4}(150,135)(140,145)
\dottedline{4}(150,137)(142,145)
\dottedline{4}(150,139)(144,145)
\dottedline{4}(150,141)(146,145)
\dottedline{4}(150,143)(148,145)

\dottedline{4}(150,145)(150,145)


\end{picture}
\vskip 1.3cm

%% file: Summary.tex
\section{Summary of Results}
\setcounter{theorem}{0}\setcounter{figure}{0}\setcounter{table}{0}


It is obvious that, for $\chi = \chizero,\ldots,\chifour$, $\P\Stab_{\Isom\Lambda}(\Fix\,\chi)$ is the subgroup of $\P\Isom(\Fix\,\chi)$ consisting of elements that extend back to isometries of the full $\Zi$-lattice $\Lambda$.  In other words, for $i = 0, \ldots, 4$,
\begin{equation*}
     \P\Stab_{\Isom\Lambda}(\Fix\,\chi_{i})
     \; \cong \;
     \P\!\left(\left\{\;M\in\Isom(L_{i})\;\left\vert\;
     \Basisi\cdot M \cdot \Basisi^{-1} \in \Zi^{6 \times 6}
     \right.\;\right\}\right)
\end{equation*}
where the $L_{i}$'s are $\Z$-lattices given as in subsection \ref{subsection:FixedLattices}, and the $\Basisi$'s in subsection \ref{subsection:FixedBases}.  Since $\P\Stab_{\Isom\Lambda}(\Fix\,\chi_{i})$ is defined by algebraic equations with coefficients in $\Zi$, it is an arithmetic subgroup of $\Isom(\RHfive)$.  Recall also that $L_{i} \cong \Fix(\chi_{i})$ as abstract $\Z$-lattices.  We may now summarize the results of this article as follows:

\begin{theorem}\mbox{}\\
\begin{enumerate}
     \item  \vskip -0.5cm
              The moduli space $\MsRi$ of stable real binary octics of type
              $i = 0, \ldots, 4$ is isomorphic as a metric space to the following
              quotient of real hyperbolic $5$-space:
              \begin{equation*}
                   \MsRi \; \cong \;
                   \P\GammaR_{i}\left\backslash\left(\RHfive_{[\chi_{i}]}\right)\right.,
              \end{equation*}     
              where $\P\GammaR_{i}$ $:=$
              $\Stab_{\P\Isom\Lambda}\left(\RHfive_{[\chi_{i}]}\right)$, and
              \begin{equation*}
                   \Stab_{\P\Isom\Lambda}\left(\RHfive_{[\chi_{i}]}\right)
                   \cong
                   \textnormal{\small$
                   \left\{\begin{array}{cl}
                        \P\Stab_{\Isom\Lambda}(\Fix\,\chi_{i}),                         & i = 0, 1, 2, 4 \\
                        \P\Stab_{\Isom\Lambda}(\Fix\,\chi_{i}) \rtimes (\Z/2\Z), & i = 3
                   \end{array}\right.$}.
              \end{equation*}
   \item   The moduli space $\MoRi$ of smooth real binary octics of type
              $i = 0, \ldots, 4$ is isomorphic as a metric space to the following
              open subspace of the above quotient:
              \begin{equation*}
                   \MoRi \; \cong \;
                   \P\GammaR_{i}\left\backslash\left(\RHfive_{[\chi_{i}]}-\H\right)\right..
              \end{equation*}     
    \item  Each\, $\P\Stab_{\Isom\Lambda}(\Fix\,\chi_{i})$, $i = 0, \ldots, 4$,
              is an arithmetic subgroup of\, $\Isom(\RHfive)$; hence, each has
              finite co-volume and is (isomorphic to) a finite-index subgroup of
              $\P\Isom(L_{i})$.  Consequently, each $\P\GammaR_{i}$ is
              commensurable with $\P\Isom(L_{i})$.
    \item  Each\; $\P\Isom(L_{i})$, $i = 0, \ldots, 4$, is a discrete reflection
              subgroup of $\Isom(\RHfive)$ with Vinberg diagram given as in
              subsection \ref{subsection:VinbergDiagrams}.
\end{enumerate}
\end{theorem}

The other main result of this paper is Corollary \ref{Corollary:MsRIsNOTHyperbolic}, which states that $\MsR$ is not a real-hyperbolic orbifold.
